\documentclass[11pt]{article}
\usepackage{amstext,amssymb,amsmath,amsbsy}

\usepackage{hyperref}
\usepackage{amscd}
\usepackage{amsfonts}
\usepackage{indentfirst}
\usepackage{verbatim}
\usepackage{amsmath}
\usepackage{amsthm}
\usepackage{cleveref}
\usepackage{enumerate}
\usepackage{graphicx}
\usepackage{subfig}
\usepackage{color}
\usepackage[OT1]{fontenc}
\usepackage[latin1]{inputenc}
\usepackage[english]{babel}
\usepackage{amssymb}
\usepackage{tikz}
\usepackage{dsfont}
\usepackage{amssymb,MnSymbol}

\usepackage{epsfig} 

\newtheorem{theorem}{Theorem}
\newtheorem{lemma}{Lemma}

\newtheorem{proposition}{Proposition}
\newtheorem{corollary}{Corollary}
\newtheorem{remark}{Remark}
\newtheorem{example}{Example}
\newtheorem{exercise}{Exercise}

\newtheorem{definition}{Definition}
\theoremstyle{definition}
\numberwithin{equation}{section}
\numberwithin{theorem}{section}
\numberwithin{lemma}{section}
\numberwithin{proposition}{section}
\numberwithin{remark}{section}
\numberwithin{example}{section}
\numberwithin{corollary}{section}
\numberwithin{exercise}{section}
\numberwithin{definition}{section}

\newcommand{\dsp}{\displaystyle}

\newcommand{\eps}{\varepsilon}

\newcommand{\mC}{\mathbb{C}}

\newcommand{\mN}{\mathbb{N}}

\newcommand{\mR}{\mathbb{R}}
\newcommand{\mP}{\mathbb{P}}

\newcommand{\hu}{\hat u}

\newcommand{\tu}{{\widetilde u}}

\newcommand{\tA}{{\widetilde A}}
\newcommand{\tB}{{\widetilde B}}
\newcommand{\tF}{{\widetilde F}}

\newcommand{\hA}{{\hat A}}

\newcommand{\hb}{{\hat b}}
\newcommand{\tb}{{\widetilde b}}

\newcommand{\cF}{{\cal F}}

\newcommand{\cG}{{\cal G}}

\newcommand{\cM}{{\cal M}}

\newcommand{\cU}{{\cal U}}
\newcommand{\cC}{{\cal C}}
\newcommand{\cL}{{\cal L}}

\newcommand{\cR}{{\cal R}}

\newcommand{\rank}{\mbox{rank }}

\newcommand{\R}{\mathbb{R}}

\newcommand{\hx}{\hat x}

\newcommand{\tvarphi}{\widetilde \varphi}

\newcommand{\tx}{\widetilde x}

\newcommand{\ty}{\widetilde y}

\newcommand{\txi}{\widetilde \xi}

\newcommand{\tdeltax}{\widetilde{\delta x}}

\newcommand{\tdeltau}{\widetilde{\delta u}}

\newcommand{\be}{\begin{equation}}
\newcommand{\ee}{\end{equation}}

\newcommand{\mH}{\mathbb{H}}

\newcommand{\mU}{\mathbb{U}}

\newcommand{\cD}{{\mathcal D}}

\title{A short introduction to the control theory in finite-dimensional spaces}

\author{Hoai-Minh Nguyen \footnote{Sorbonne Universit\'e, Universits\'e Paris Cit\'e, CNRS, INRIA, Laboratoire Jacques-Louis Lions, LJLL, F-75005 Paris, France, hoai-minh.nguyen@sorbonne-universite.fr.}}

\begin{document}

\maketitle 

\tableofcontents

\section{Introduction} This is a brief introduction to control theory in finite-dimensional spaces. The material is partly based on my lectures for the Master 1 program in Math\'ematiques et applications at Sorbonne University, delivered over the past few years. 
The aim is to provide a concise overview of the subject, primarily focusing on the linear setting. Proofs are presented in detail and are selected to allow for extensions to the infinite-dimensional case in many situations. Topics covered include the Kalman rank condition, the Hautus test, observability, stability, 
detectability and dynamic observers,  the Pole Shifting Theorem, the linear test for controllability, linear-quadratic optimal control over finite and infinite horizons, and stabilization via Gramians. The references for these materials and much more can be found in \cite{lee1967foundations,kwakernaak1972linear,Sontag98,Coron07,Zabczyk08,brockett2015finite,trelat2024control}. 
The author thanks Professor Jean-Michel Coron for introducing him to control theory, as well as for many insightful discussions and valuable guidance. 

One of the motivations for presenting the material from my classes in the form of lecture notes comes from my talk at the conference ``Resonances in the Mathematical World", held in the summer of 2024 at the University of Science, Vietnam National University - Ho Chi Minh City. The event was a special occasion to celebrate the 60th birthday of Professor Dang Duc Trong, a respected teacher of generations of students from the Faculty of Mathematics and Computer Science. The goal of my talk was to briefly present my recent work on stabilization in \cite{Ng-Riccati, Ng-S-KdV, Ng-S-Schrodinger} with the roots from the optimal control theory. While this lecture does not reproduce the content of those works, it does highlight some of the ideas. Readers interested in the details may consult those papers after becoming familiar with the theory of strongly continuous semigroups. The author also thanks the organizers of the conference for hosting such a meaningful event and the editors of this special issue for giving him the opportunity to present the lecture. These notes are dedicated to Professor Dang Duc Trong.

\section{Controllability - the Kalman rank condition}

Let $n \ge 1$ and $p \ge 1$. We first deal with autonomous setting. Consider the control system 
\begin{equation}\label{sys-LA}
x' = A x + Bu, 
\end{equation}
where $x \in \mR^n$, $u \in \mR^p$, and $A \in \mR^{n \times n}$ and $B \in \mR^{n \times p}$ are two constant matrices. 

\medskip
Here is the basic control question for this system: Let $t_0 < t_1$. Given $x_0, x_1 \in \mR^n$, can one find a control $u$ (in an appropriate space) such that the unique solution of \eqref{sys-LA} with $x(t_0) = x_0$ satisfying $x(t_1) = x_1$? 

\medskip 
Before answering this question, let us state and prove the following standard well-posedness results related to \eqref{sys-LA} in which the matrix $A$ might depend on $t$. 

\begin{lemma}\label{lem-CL1} Let $t_0 < t_1$, $A \in L^\infty((t_0, t_1); \mR^{n \times n})$, and $b \in L^1((t_0, t_1); \mR^n)$. Given $x_0 \in \mR^n$, there exists a unique solution $x \in C([t_0, t_1]; \mR^n)$ of the system 
\be \label{lem-CL1-sys}
x'(t) = A(t) x(t) + b(t) \mbox{ in } [t_0, t_1] \quad \mbox{ and } \quad x(t_0) = x_0, 
\ee
in the sense that 
\begin{equation}\label{lem-CL1-sol}
x(t) = x_0 + \int_{t_0}^t A(s) x(s) \, ds + \int_{t_0}^t b(s) \, ds \mbox{ for } t \in [t_0, t_1]. 
\end{equation}
\end{lemma}

\begin{remark} \rm Identity \eqref{lem-CL1-sol} can be viewed as a weak form of \eqref{lem-CL1-sys}. In fact, if $A$ and $b$ are smooth and $x$ is a smooth solution of \eqref{lem-CL1-sys}, then \eqref{lem-CL1-sol} holds after integrating the equation of $x$ with respect to time from $t_0$ to $t$. 
\end{remark}

\begin{proof} We first give the proof in the case $t_1 - t_0 $ is sufficiently small. We then show how to modify to obtain the proof in the general case. 

Assume that $t_1 - t_0 $ is sufficiently small. Consider the map  
$$
\begin{array}{cccc}
\cF: & C([t_0, t_1]; \mR^n) & \mapsto & C([t_0, t_1]; \mR^n) \\[6pt]
& x & \to & \cF x 
\end{array}
$$
defined by 
$$
\cF x(t): = x_0 + \int_{t_0}^t A(s) x(s) \, ds + \int_{t_0}^t b(s) \, ds \mbox{ for } t \in [t_0, t_1]. 
$$
In this part of the proof, $C([t_0, t_1]; \mR^n)$ is equipped with the standard sup norm: 
$$
\| f\|_{L^\infty(t_0, t_1)}: = \sup_{t \in [t_0, t_1]} \|f(s)\| \mbox{ for } f \in C([t_0, t_1]; \mR^n), 
$$
where $\|f(s)\|$ denotes the Euclidean norm in $\mR^n$ of $f(s)$.

We are proving that $\cF$ is a contraction mapping, which ensures the existence and uniqueness of the solution, as the space under consideration is a Banach space.  We have, by the definition of $\cF$, for $x, y \in C([t_0, t_1]; \mR^n)$, 
\begin{equation*}
\cF x (t) - \cF y (t) =  \int_{t_0}^t A(s) \big( x(s) - y(s) \big) \, ds. 
\end{equation*}
This yields 
\begin{equation} \label{lem-CL1-**}
\| \cF x (t) - \cF y (t) \| \le \int_{t_0}^t \| A(s) \| \|x(s) - y(s) \| \, ds \le \int_{t_0}^t \| A \|_{L^\infty} \|x - y \|_{L^\infty} \, ds. 
\end{equation}

Here and in what follows, one sometimes just denotes $\| \cdot \|$ the norm of a quantity in an appropriate space  for notational ease. For example, in the above inequality, $\| \cF x (t)  - \cF y (t) \|$ is the Euclidean norm of 
$ \cF x (t) - \cF y (t)$, $\| A(s) \|$ is the matrix norm of $A(s)$ \footnote{Here, we considered $A(s)$ as a linear application from $\mR^n$ to $\R^n$ equipped the Euclidean norm, and $\| A(s) \|$ denotes the corresponding norm. Other norms in $\mR^n$ or other norms of $A(s)$ can be also considered and the inequality holds with an additional positive constant on the RHS. This fact does not change the process of the proof.},  
$\|x(s) - y(s) \|$ is the Euclidean norm of  $x (s) - y (s)$. 

It follows from \eqref{lem-CL1-**} that 
\begin{equation*}
\| \cF x  - \cF y \|_{L^\infty(t_0, t_1)} \le (t_1- t_0) \| A \|_{L^\infty(t_0, t_1)} \|x - y \|_{L^\infty (t_0, t_1)}, 
\end{equation*}
so $\cF$ is contracting if $t_1 - t_0$ is sufficiently small.

To deal with the general case, one can use two methods. The first one is to extend the time interval consecutively. The details for this are omitted. The second one is to modify the norm in $C([t_0, t_1]; \mR^n)$ as follows. Consider the norm, for large positive $\lambda$,  
$$
\| x  \|_{L_\lambda^\infty(t_0, t_1)} : = \sup_{t \in [t_0, t_1]} e^{-\lambda t} \| x(t) \|. 
$$
One then has, for $t \in [t_0, t_1]$,  
\begin{multline*}
\| \cF x (t) - \cF y (t) \| \le \int_{t_0}^t \| A(s) \| \|x(s) - y(s) \| \, ds \\[6pt] \le \int_{t_0}^t e^{\lambda t} \| A \|_{L^\infty(t_0, t_1)} \|x - y  \|_{L_\lambda^\infty(t_0, t_1)} \, ds 
\le \frac{e^{\lambda t}}{\lambda} \| A \|_{L^\infty(t_0, t_1)} \|x - y  \|_{L_\lambda^\infty(t_0, t_1)}, 
\end{multline*}
which yields, for $t \in [t_0, t_1]$,  
$$
\| \cF x (t) - \cF y (t) \| e^{-\lambda t} \le \frac{1}{\lambda} \| A \|_{L^\infty(t_0, t_1)} \|x - y  \|_{L_\lambda^\infty (t_0, t_1)}. 
$$
In other words, 
$$
\| \cF x  - \cF y  \|_{L^\infty_\lambda(t_0, t_1)} \le \frac{1}{\lambda} \| A \|_{L^\infty(t_0, t_1)} \|x - y  \|_{L_\lambda^\infty(t_0, t_1)}. 
$$
Thus $\cF$ is contracting for large positive $\lambda$. The conclusion follows. 
\end{proof}

We now go back to the control question. Concerning this question, one has the following definition in which the functional space of $u$ is given.  

\begin{definition} \rm  Let $t_0 < t_1$. System \eqref{sys-LA} is controllable for the time interval $[t_0, t_1]$ if for all $x_0, x_1 \in \mR^n$, there exists $u \in L^2((t_0, t_1), \mR^p)$ such that the solution of \eqref{sys-LA} with $x(t_0) = x_0$ satisfying $x(t_1) = x_1$. 
\end{definition}


Some conditions on $(A, B)$ are required for the controllability of \eqref{sys-LA}. Indeed, if $A = 0$ and $B=0$ then $x(t_1) = x_0$ for all $u$. Another less trivial example but in the same spirit is that if the last row of $A$ and $B$ are 0, then the last component of $x(t_1)$ is equal to the last component of $x(t_0)$ and the controllability is not achieved. 


\begin{remark} \rm As seen from the proof of \Cref{thm-LA-Kalman} below, the condition on the controllability of \eqref{sys-LA} is unchanged if instead of requiring $u \in L^2((t_0, t_1), \mR^p)$, one uses $u \in L^\infty((t_0, t_1), \mR^p)$ or $u \in L^1((t_0, t_1), \mR^p)$.
\end{remark}

Here is the answer to the controllability of \eqref{sys-LA} due to Kalman \cite{kalman1960contributions} and surprisingly appeared only in 1960.

\begin{theorem}[Kalman's theorem] \label{thm-LA-Kalman} System \eqref{sys-LA} is controllable if and only if $\rank M = n$, where 
\be\label{def-M}
M := [B, AB, \dots, A^{n-1} B]. 
\ee
\end{theorem}

\begin{remark} \rm Note that $M$ is an $n \times np$ matrix. \end{remark}

\begin{remark} \rm In what follows, we often mention that the pair $(A, B)$ is controllable instead of saying that System \eqref{sys-LA} is controllable for the simplicity. 
\end{remark}

Before giving the proof of \Cref{thm-LA-Kalman}, we recall that the solution of \eqref{sys-LA} with $x(t_0 ) = x_0$ is given by 
\begin{equation}\label{sys-LA-sol}
x(t) = e^{(t-t_0) A} x_0 + \int_{t_0}^t e^{(t-s) A} B u(s) \, ds. 
\end{equation}

Given a real matrix $C$, we will denote $C^\top$ its transpose. 

\begin{proof}[Proof of \Cref{thm-LA-Kalman}]
We first motivate the way to obtain the Kalman rank condition. The starting key idea is the following: the information of a vector $\xi$ in $\mR^n$ is equivalent to the information of all $\langle \xi, y \rangle $ for $y \in \mR^n$, where $\langle \cdot, \cdot \rangle$ denotes the Euclidean scalar product \footnote{This idea is also the basis of the Hilbert uniqueness method applied also to the infinite dimensional setting, see, e.g., \cite{Lions88-VolI}.}. 

From \eqref{sys-LA-sol}, we have, for $y \in \mR^n$,  
\begin{equation*}
\langle x(t_1) - e^{(t_1 - t_0) A} x_0, y \rangle = \int_{t_0}^{t_1} \langle e^{(t_1-s) A} B u(s), y \rangle  \, ds,  
\end{equation*}
which yields  
\begin{equation} \label{thm-LA-Kalman-p1}
\langle x(t_1) - e^{(t_1 - t_0) A} x_0, y \rangle = \int_{t_0}^{t_1} \langle u(s), B^\top e^{(t_1-s) A^\top}y \rangle  \, ds. 
\end{equation}

Denote
\begin{equation}\label{thm-LA-Kalman-def-cM}
\cM: = \Big\{g \in L^2((t_0, t_1); \mR^p) ; \;  g(s) =  B^\top e^{(t_1-s) A^\top} y \mbox{ in } [t_0, t_1] \mbox{ for some } y \in \mR^n\Big\},
\end{equation}
and equip this space with the standard $L^2((t_0, t_1); \mR^p) $-scalar product. Note that $\dim \cM \le n $ so $\cM $ is closed. Denote $\mP_{\cM}$ the projection from $L^2((t_0, t_1); \mR^p)$ into $\cM$. 

Another simple but important ingredient of the proof is the following identity, which is a consequence of \eqref{thm-LA-Kalman-p1} and the property of the projection, 
\begin{equation} \label{thm-LA-Kalman-p2}
\langle x(t_1) - e^{(t_1 - t_0) A} x_0, y \rangle = \int_{t_0}^{t_1} \langle \mP_{\cM} u (s), B^\top e^{(t_1-s) A^\top}y \rangle  \, ds. 
\end{equation}

Let $z \in \mR^n$ be such that 
$$
\mP_{\cM} u (s) = B^\top e^{(t_1-s) A^\top} z \mbox{ in } (t_0, t_1). 
$$
It follows from \eqref{thm-LA-Kalman-p2} that 
\begin{equation*}
\langle x(t_1) - e^{(t_1 - t_0) A} x_0, y \rangle = \int_{t_0}^{t_1} \langle B^\top e^{(t_1-s) A^\top} z, B^\top e^{(t_1-s) A^\top} y \rangle  \, ds. 
\end{equation*}
We thus obtain, for $y \in \mR^n$,  
\begin{equation} \label{thm-LA-Kalman-p2-**}
\langle x(t_1) - e^{(t_1 - t_0) A} x_0, y \rangle = \int_{t_0}^{t_1} \langle e^{(t_1 - s) A} BB^\top e^{(t_1-s) A^\top} z,  y \rangle  \, ds. 
\end{equation}
 
Set  
 \begin{equation}\label{thm-LA-Kalman-def-mC}
\mC : = \int_{t_0}^{t_1} e^{(t_1-s) A} B B^\top e^{(t_1 -s) A^\top} \, ds. \end{equation}
It is clear that $\mC$ is symmetric and $\mC$ is non-negative; moreover, 
$$
\langle \mC z, z \rangle = \int_{t_0}^{t_1} \|B^\top e^{(t_1 -s) A^\top} z\|^2 \, ds.
$$

Identity \eqref{thm-LA-Kalman-p2-**} can be written as
\be
\label{thm-LA-Kalman-p3}
\langle x(t_1) - e^{(t_1 - t_0) A} x_0, y \rangle = \langle \mC z, y \rangle \mbox{ for all } y \in \mR^n.  
\ee

We have thus proved that for $u \in L^2((t_0, t_1); \mR^n)$, it holds
\be
\label{thm-LA-Kalman-p4}
x(t_1) - e^{(t_1 - t_0) A} x_0  =  \mC z, 
\ee 
where $x$ is the solution of \eqref{sys-LA} with $x(t_0) = x_0$, and $z \in \mR^n$ is such that 
$$
\mP_{\cM} u (s) = B^\top e^{(t_1-s) A^\top} z \mbox{ in } (t_0, t_1). 
$$

Identity \eqref{thm-LA-Kalman-p4} is the key point in the proof given below to derive the Kalman rank condition. It yields the following important fact
\be \label{thm-LA-Kalman-p5}
\mbox{system \eqref{sys-LA} is controllable if and only if the matrix $\mC$ is invertible.} 
\ee

Indeed, if $\mC$ is invertible, then for all $x_1 \in \mR^n$, let $u$ be defined by 
\be \label{thm-LA-Kalman-def-u}
u(s) = B^\top e^{(t_1- s) A^\top} z \mbox{ for } s \in (t_0, t_1), \mbox{ where } z = \mC^{-1} (x_1 - e^{(t_1 - t_0) A} x_0).  
\ee
By \eqref{thm-LA-Kalman-p4}, $u$ is a control steering $x_0$ at time $t_0$ to $x_1$ at time $t_1$. 

We next show that if the system is controllable, then $\mC$ is invertible by contradiction. Assume that $\mC$ is not invertible. Then by taking $x_0 = 0$ and $x_1 \in (\mbox{range } \mC)^\perp \setminus \{0 \}$ \footnote{Here $\mbox{range } \mC : = \{\mC z; z \in \mR^n \} \subset \mR^n$ and the perpendicular relation is considered with respect to the Euclidean scalar product.}, one derives that there is no $z \in \mR^n$ such that \eqref{thm-LA-Kalman-p4} holds. Thus \eqref{thm-LA-Kalman-p5} is proved.

We next show that the invertibility of $\mC$ is equivalent to the fact $\rank M = n$, and the conclusion follows.

Since $\mC$ is symmetric and $\mC \ge 0$, it follows that 
\be \label{thm-LA-Kalman-c1}
\mbox{ the invertibility of $\mC$ is equivalent to the fact $\mC > 0$},  
\ee
i.e., \footnote{This inequality is known as the observability inequality of the control system $x' = A x + B u$ in the time interval $[t_0, t_1]$.}
\be \label{thm-LA-Kalman-obs}
\langle \mC z, z \rangle = \int_{t_0}^{t_1} \|B^\top e^{(t_1 -s) A^\top} z\|^2 \, ds \ge C_{t_0, t_1} \| z\|^2 \mbox{ for all } z \in \mR^n, 
\ee
for some positive constant $C_{t_0, t_1}$.

\medskip
The proof is now divided into two steps. 

\medskip 

\noindent {\it Step 1:} The invertibility of $\mC$ implies the rank condition.   

We prove this by contradiction. Assume that $\rank M < n$ and thus $\rank M^\top < n$.  There exists $y \in \mR^n \setminus \{ 0 \}$ such that 
\be \label{thm-LA-Kalman-p6}
B^\top (A^\top)^k y = 0 \mbox{ for all } 0 \le k \le n-1.  
\ee
It follows from the Caley-Hamilton theorem that  $(A^{T})^k$ for $k \ge 0$ is a linear combination of $I, A^\top, \dots, (A^\top)^{n-1}$. We derive from \eqref{thm-LA-Kalman-p6} that 
\be \label{thm-LA-Kalman-p7}
B^\top (A^\top)^k y = 0 \mbox{ for all } k \ge 0.  
\ee

Set 
\be \label{thm-LA-Kalman-def-f}
f(t) = B^\top e^{t A^\top} y \mbox{ for } t \in \mR. 
\ee
Then $f$ is analytic on $\mR$ and 
$$
f^{k} (0) = B^\top (A^\top)^k y \mbox{ for } k \ge 0. 
$$
We derive from \eqref{thm-LA-Kalman-p7} and the analyticity of $f$ that 
$$
f(t) = 0 \mbox{ for }  t \in \mR.   
$$
This implies 
$$
\langle \mC y, y \rangle =  \int_{t_0}^{t_1} \|B^\top e^{(t_1 -s) A^\top} y\|^2 \, ds =   \int_{t_0}^{t_1}  \| f(t_1 - s) \|^2 \, ds = 0. 
$$
This contradicts the invertibility of $\mC$ by \eqref{thm-LA-Kalman-c1}. 

\medskip 
\noindent{\it Step 2:} The rank condition implies the invertibility of $\mC$. 

For $y \in \mR^n \setminus \{0\}$. Since $\rank M^\top = n$, it follows that 
$$
B^\top (A^\top)^k y \neq 0 \mbox{ for some } 0 \le k \le n-1. 
$$
Thus the function $f$ defined by \eqref{thm-LA-Kalman-def-f} is not identically equal to $0$ by the analyticity of $f$. Using again the analyticity of $f$, we derive that 
$$
\langle \mC y, y \rangle =  \int_{t_0}^{t_1} \|B^\top e^{(t_1 -s) A^\top} y\|^2 \, ds =   \int_{t_0}^{t_1}  \| f(t_1 - s) \|^2 \, ds >0. 
$$
Thus $C$ is invertible. 

The proof is complete. 
\end{proof}

\begin{remark} \rm The proof given above is in the spirit of the standard proof of the Kalman rank condition. Nevertheless, it is presented in such a way that the Kalman rank condition naturally appears via the invertibility of a matrix $\mC$ given in \eqref{thm-LA-Kalman-def-mC}. This invertibility is equivalent to an observability inequality \eqref{thm-LA-Kalman-obs}. 
\end{remark}

It is convenient to summary the result given in \eqref{thm-LA-Kalman-p4}. Let $t_0 < t_1$ and set 
\begin{equation}\label{lem-Kalman-def-mC}
\mC  = \int_{t_0}^{t_1} e^{(t_1-s) A} B B^\top e^{(t_1 -s) A^\top} \, ds. 
\end{equation}
We have 

\begin{lemma}\label{lem-Kalman}  Let $t_0 < t_1$ and $x_0 \in \mR^n$. The following two facts hold: 

\begin{itemize}
\item[i)] Given $u \in L^2((t_0, t_1); \mR^n)$, let  $x$ be the solution of the system 
\be \label{lem-Kalman-sys}
x'  = Ax + Bu \mbox{ in } (t_0, t_1) \quad \mbox{ and } \quad x(t_0) = x_0. 
\ee
There exists $z \in \mR^n$ such that
\be \label{lem-Kalman-identity}
x(t_1) - e^{(t_1 - t_0) A} x_0  =  \mC z. 
\ee 

\item[ii)] Given $z \in \mR^n$, and set 
$$
u(s) = B^\top e^{(t_1-s)A^\top} z \mbox{ in } (t_0, t_1). 
$$
Then \eqref{lem-Kalman-identity} holds where $x$ is the solution of \eqref{lem-Kalman-sys}.
\end{itemize}
\end{lemma}

The given proof also reveals the following interesting property of the control chosen in the proof of \Cref{thm-LA-Kalman} in \eqref{thm-LA-Kalman-def-u}. 

\begin{proposition} Let $t_0 < t_1$ and  $x_0, x_1 \in \mR^n$. Assume that $\rank M = n$ and set 
\begin{equation*}
\tu(t) = B^\top e^{(t_1-t) A^\top} z \quad \mbox{ with } \quad  z = \mC^{-1} (x_1 - e^{(t_1 - t_0) A} x_0), 
\end{equation*}
where $\mC$ is defined in \eqref{lem-Kalman-def-mC}.  Then $\tu$ steers the system from $x_0$ at time $t_0$ to $x_1$ at time $t_1$. Moreover, if $u \in  L^2((t_0, t_1); \mR^p)$ is any control which steers the system from $x_0$ at time $t_0$ to $x_1$ at time $t_1$, then 
\begin{equation*}
\int_{t_0}^{t_1} |\tu(t)|^2 \, dt \le \int_{t_0}^{t_1} |u(t)|^2 \, dt. 
\end{equation*}
The equality holds if and only if $u = \tu$ a.e. in $[t_0, t_1]$. 
\end{proposition}

\begin{example}[Pendulum example] \rm The pendulum can be model by 
\begin{equation*}
m \theta'' + m g \sin \theta = u,  
\end{equation*}
where $\theta$ is the angle with respect to the horizontal line and $u$ is the torque. There are two equilibrium positions for this system. One is $(\theta, u) = (0, 0)$ which is stable and one is $(\theta, u) = (\pi, 0)$ which is unstable. 
We are now interested in states for which $\theta$ is around $\pi$. Set $\varphi = \theta - \pi$ and assume that all constants $m$ and $g$ are 1 for notational ease. One has 
\begin{equation*}
\varphi'' - \sin \varphi = u. 
\end{equation*}
The linearized system around $(0, 0)$ is 
\be\label{pendulum-pi}
\varphi'' - \varphi = u. 
\ee
This system can be written under the form 
$$
x' = A x + B u \quad \mbox{ where }  \quad  x = \left(\begin{array}{cc}
\varphi \\
\varphi'
\end{array}\right), \quad  A = \left(\begin{array}{cc} 0 & 1 \\
1 & 0
\end{array}\right), \quad B = \left(\begin{array}{cc} 0\\
1 
\end{array}\right). 
$$
One can check that the linearized system is controllable since 
$$
M = [B, AB]= \left(\begin{array}{cc} 0 & 1 \\
1 & 0
\end{array}\right), \mbox{ whose rank is $2$.}
$$
\end{example}

It is natural to ask what happens in the case $\rank M < n$.  Here is a result in this direction dealing with $x_0 = 0$. 

\begin{proposition}\label{pro-C1} Let $t_0 < t_1$ and set 
$$
\cC = \mbox{range } \mC = \{\mC x ; \, x \in \mR^n \} \quad \mbox{ and } \quad  \cU = \cC^\perp. 
$$
Then 
\begin{itemize}
\item[i)] 
$$
\mR^n = \cC \otimes \cU. 
$$

\item[ii)] One can steer $0$ at time $t_0$ to every element in $\cC$ at time $t_1$. 

\item[iii)] One cannot steer $0$ at time $t_0$ to any element in $\cU \setminus \{0 \}$ at time $t_1$. 

\end{itemize}
Moreover, the space $\cU$ can be characterized by 
\be \label{pro-C1-cl1}
\cU : =  \mbox{ ker } M^\top = \Big\{ y \in \mR^n; \; M^\top y = 0 \Big\},   
\ee
where $M$ is defined in \eqref{def-M}. 
\end{proposition}

\begin{proof} Assertion $i)$, $ii)$, and $iii)$ follows from \Cref{lem-Kalman}. It remains to prove \eqref{pro-C1-cl1}. One has
$$
\cU = \mbox{ ker } \mC^\top = \mbox{ ker } \mC,   
$$
since $\mC$ is symmetric.  One notes that 
$$
\mbox{ ker } \mC =  \Big\{ y \in \mR^n; \; M^\top y = 0 \Big\}. 
$$
It follows that 
$$
\cU = \mbox{ ker } M^\top. 
$$
The proof is complete. 
\end{proof}


\begin{lemma} \label{lem-C1} Let $b_j$ denote the $j$-th column of $B$ and define 
$$
\cR(A, B) : = \mbox{span } \Big\{ A^i b_j;  \; i \ge 0, 1 \le j \le p \Big\}, 
$$
which is called the controllable space. Then 
$$
\cC = \cR(A, B). 
$$
\end{lemma}

\begin{remark} \label{rem-AB-invariant} \rm $\cR(A, B)$ is the smallest $A$-invariant subspace of $\mR^n$ containing $B$. 
\end{remark}

\begin{proof} From \Cref{pro-C1}, one has 
$$
\cC = \Big(\mbox{ ker } M^\top  \Big)^\perp = \Big\{ y \in \mR^n; \; M^\top y = 0 \Big\}^\perp. 
$$
It follows that  $
\cC = \mbox{range } M$ and the conclusion follows. 
\end{proof}

We also have the following important result. 
 
\begin{lemma}[Kalman-controllability condition] \label{lem-Kalman-cond} Assume that $(A, B)$ is not controllable. Denote $r = \dim \cR(A, B) < n$. There exists an invertible $n \times n$ matrix $T$ such that, for 
\be \label{lem-Kalman-cond-tAtB}
\widetilde A : = T^{-1} A  T \quad \mbox{ and } \quad \widetilde B : = T^{-1} B, 
\ee
we have the following block structures
\begin{equation}\label{lem-Kalman-cond-structure}
\widetilde A  =  \left(\begin{array}{cc} A_1 & A_2 \\
0 & A_3
\end{array}\right),  \quad \widetilde B  =  \left(\begin{array}{c} B_1  \\
0 
\end{array}\right), 
\end{equation}
with $A_1 \in \mR^{r \times r}$ and $B_1 \times \mR^{r \times m}$. Moreover, $(A_1, B_1)$ is controllable.  
\end{lemma}

\begin{remark} \rm Some comments on \Cref{lem-Kalman-cond} are in order. Let $x$ be a solution of the system 
$$
x' = \tA x + \tB u
$$
and denote $x = (x_1, x_2)^\top \in \mR^{r} \times \mR^{n-r}$. One then has, by \eqref{lem-Kalman-cond-structure},  
$$
\left\{\begin{array}{c}
x_1' = A_1 x_1 + A_2 x_2 + B_1 u, \\[6pt]
x_2' = A_3 x_2.  
\end{array}\right. 
$$
Hence,  $x_1 = y_1 + z_1$, where 
$$
\left\{\begin{array}{c}
z_1' = A_2 x_2, \quad z_1(t_0) = 0, \\[6pt]
y_1' = A_1 y_1 + B_1 u, \quad y_1(t_0) = x_1(t_0). 
\end{array}\right. 
$$
\Cref{lem-Kalman-cond} says that, after a change of variables, the first $r$ components $x_1$  of $x$ is controllable and the last $n-r$ components $x_2$ of $x$ is uncontrollable.  
\end{remark}

\begin{proof} We have, by \Cref{pro-C1} and \Cref{lem-C1},  
$$
\mR^n = \cR \otimes \cU, 
$$
where $\cR = \cR(A, B)$ and $\cU = \cR^\perp$. Let $v_1, \cdots, v_r$ and $v_{r+1}, \cdots, v_{n}$ be an orthogonal basis of $\cR$ and $\cU$, respectively. Then $\{v_1, \cdots, v_n\}$ is an orthogonal basis of $\mR^n$. 

Set 
\be \label{lem-Kalman-cond-def-T}
T = (v_1, \cdots, v_r, v_{r+1}, \cdots, v_{n}) \in \mR^{n \times n}. 
\ee
Then  
\be \label{lem-Kalman-cond-T}
T^\top = T^{-1}. 
\ee
By \eqref{lem-Kalman-cond-tAtB}, we have, for $1 \le i \le r < r+1 \le j  \le n$,  
$$
\langle e_j, \widetilde A e_i \rangle = \langle e_j, T^{-1} A T e_i \rangle = \langle e_j, T^{-1} A v_i  \rangle  \mathop{=}^{\eqref{lem-Kalman-cond-T}} \langle T e_j, A v_i  \rangle
\mathop{=}^{\eqref{lem-Kalman-cond-def-T}}\langle v_j, A v_i  \rangle  \mathop{=}^{\Cref{rem-AB-invariant}} 0,
$$
and, for $u \in \mR^p$,  
$$
\langle \widetilde B u, e_j \rangle = \langle T^{-1} B u, e_j \rangle  \mathop{=}^{\eqref{lem-Kalman-cond-T}}   \langle  B u, T e_j \rangle \mathop{=}^{\eqref{lem-Kalman-cond-def-T}}  \langle B u, v_j \rangle = 0. 
$$
Here in the last identity, we used the fact $Bu \in \cR$, $v_j \in \cU$ for $r+1 \le j \le n$, and $\cR \perp \cU$. 
This implies the structures of $\widetilde A$ and $\widetilde B$ in \eqref{lem-Kalman-cond-structure}. 

It remains to prove that $(A_1, B_1)$ is controllable. This follows from the fact 
$$
M = T \left(\begin{array}{c} M_1 \\
0 \end{array}\right), 
$$
where $M$ is defined by \eqref{def-M} and $M_1  = [B_1, A_1 B_1, \dots, A_1^{n-1} B_1]$ since 
$$
T^{-1} A^{k} B = (T^{-1}A T)^k T^{-1} B  = \tA^k \tB.
$$ 
The proof is complete. 
\end{proof}

\begin{exercise} Assume that  $m=1$, $A = \mbox{diag }(\lambda_1, \cdots, \lambda_n)$,  and $B = (b_1, \cdots, b_n)^\top$. Then $(A, B)$ is controllable if and only if 
$$
\lambda_i \neq \lambda_j \mbox{ for } i \neq j \quad \mbox{ and } \quad b_i \neq 0 \mbox{ for all $i$}. 
$$
\end{exercise}

We next discuss briefly the non-autonomous setting: 
\begin{equation}\label{sys-LNA}
x'(t) = A (t) x(t) + B (t) u(t) \mbox{ in } (t_0, t_1), 
\end{equation}
where $x \in \mR^n$, $u \in \mR^p$, and $A \in L^\infty((t_0, t_1); \mR^{n \times n})$ and $B \in L^\infty((t_0, t_1); \mR^{n \times p})$. 

\medskip 
We are ready to formula and answer the question on the controllability. 

\begin{definition} \rm  Let $t_0 < t_1$. System \eqref{sys-LNA} is controllable for the time interval $[t_0, t_1]$ if for all $x_0, x_1 \in \mR^n$, there exists $u \in L^2((t_0, t_1), \mR^p)$ such that the solution of \eqref{sys-LA} with $x(t_0) = x_0$ satisfying $x(t_1) = x_1$. 
\end{definition}

Before giving the answer, let us introduce the matrix 
\begin{equation} \label{def-mC-NA}
\mC: = \int_{t_0}^{t_1} R(t_1, s) B(s) B(s)^\top R(t_1, s)^\top \, ds, 
\end{equation}
where $R \in C([t_0, t_1]^2; \mR^{n \times n}))$ is the resolvent of the system $x'(t) = A (t) x(t)$. This is the counter part of \eqref{lem-Kalman-def-mC} in the non-autonomous setting.

Recall that for $t_0 \le s, t \le t_1$, 
$$
\partial_t R(t, s) = A (t) R(t, s) \quad \mbox{ and } \quad R(s, s) = I \mbox{ which is the identity matrix in $\mR^{n \times n}$}.  
$$
It is known that the solution of \eqref{sys-LNA} with $x(t_0) = x_0 \in \mR^n$ is unique and given by 
$$
x(t) = R(t, t_0) x_0 + \int_{t_0}^t R(t, s) B(s) u (s) \, ds. 
$$
This formula is known as Duhamel's principle. 

We are ready to state and give the proof of the controllability in the non-autonomous setting, due to Kalman again \cite{kalman1960contributions}.  

\begin{theorem}[Kalman's theorem] \label{thm-LNA-Kalman} System \eqref{sys-LNA} is controllable in the time interval $[t_0, t_1]$ if and only if the matrix $\mC$ defined in 
\eqref{def-mC-NA} is invertible.
\end{theorem}

\Cref{thm-LNA-Kalman} is a consequence of $i)$ of \Cref{lem-Kalman-NA}  below, whose proof is similar to the one of \Cref{lem-Kalman} and is omitted. 

\begin{lemma}\label{lem-Kalman-NA}  Let $t_0 < t_1$ and $x_0 \in \mR^n$. The following two facts hold: 

\begin{itemize}
\item[i)] Given $u \in L^2((t_0, t_1); \mR^n)$, let  $x$ be the solution of the system 
\be \label{lem-Kalman-sys-NA}
x(t)'  = A(t)x + B(t)u \mbox{ in } (t_0, t_1) \quad \mbox{ and } \quad x(t_0) = x_0. 
\ee
There exists $z \in \mR^n$ such that
\be \label{lem-Kalman-identity-NA}
x(t_1) - R(t_1, t_0) x_0  =  \mC z, 
\ee 
where $\mC$ is defined by \eqref{def-mC-NA}. 

\item[ii)] Given $z \in \mR^n$, and set 
$$
u(s) = B(s)^\top  R(t_1, s)^\top z  \mbox{ in } (t_0, t_1). 
$$
Then \eqref{lem-Kalman-identity-NA} holds where $x$ is the solution of \eqref{lem-Kalman-sys-NA}
\end{itemize}
\end{lemma}

We also have 

\begin{proposition} Let $t_0 < t_1$ and  $x_0, x_1 \in \mR^n$. Assume that System \eqref{sys-LNA} is controllable for the time interval $[t_0, t_1]$ and set 
\begin{equation*}
\tu(t) = B(s)^\top R(t_1, s)^\top y  \quad \mbox{ with } \quad  z = \mC^{-1} (x_1 - R(t_1, t_0) x_0). 
\end{equation*}
Then the control $\tu$ steers the system from $x_0$ at time $t_0$ to $x_1$ at time $t_1$. Moreover, if $u \in L^2(t_0, t_1); \mR^p)$  is any control which steers the system from $x_0$ at time $t_0$ to $x_1$ at time $t_1$, then 
\begin{equation*}
\int_{t_0}^{t_1} |\tu(t)|^2 \, dt \le \int_{t_0}^{t_1} |u(t)|^2 \, dt. 
\end{equation*}
The equality holds if and only if $u = \tu$ a.e. in $[t_0, t_1]$. 
\end{proposition}

Here is an interesting consequence of \Cref{thm-LNA-Kalman}. 

\begin{corollary} \label{cor-LNA-Kalman} Let $t_0 < t_1$ and assume that the control system is analytic, i.e., $A$ and $B$ are analytic in an open interval containing $[t_0, t_1]$.  Assume that the control system is controllable for the time interval $[t_0, t_1]$. Then the control system is controllable for every interval $[t_0, t_0 + \delta] \subset [t_0, t_1]$. 
\end{corollary}

\begin{remark} \rm
The necessary and sufficient condition for controllability given in \Cref{thm-LNA-Kalman} requires computing the matrix $\mC$, which might be quite difficult (and even
impossible) in many cases. There are other criterion which overcome this, see, e.g., \cite[Section 1.3]{Coron07}.  
\end{remark}

\begin{remark} \rm  The finite-dimensional setting \eqref{sys-LA} can be extended to infinite-dimensional one to deal with partial differential equations. In this context, 
the state $x$ belongs to a Hilbert space $\mH$, the control $u$ belongs to a Hilbert space $\mU$, and  $A$ is a generated of a strongly continuous defined in $\mH$ (which is unbounded in general with the domain $\cD(A)$) and $B$ is an linear application from $\mU$ to a space containing $\cD(A^*)'$ (the dual space of $\cD(A^*)$) where $A^*$ is the adjoint of $A$. Thus $B$ can be bounded from $\mU$ to $\mH$ or not. There are now various notions of controllability such as null-controllability, exact controllability, and approximative controllability. These various controllability properties are then proved to be equivalent to some observability inequalities. In general, these notions are not equivalent and this is contrast to the finite dimensional case. 
The materials in this directions can be found in, see, e.g., \cite{CZ95,Coron07,BDDM07,Zabczyk08,TW09}. Tools from functional analysis, see, e.g., \cite{Brezis-FA} are indispensable in this direction. 
\end{remark}

\begin{remark} \rm The reader interested in similar phenomena as in \Cref{cor-LNA-Kalman} in the infinite dimensional spaces is referred to \cite{coron2021optimal} dealing with the controllability of hyperbolic systems in one dimensional space.
\end{remark}

\section{The Hautus test}

We next discuss the Hautus test \cite{hautus1969controllability} which gives another way to check the controllability of the pair $(A, B)$.  

\begin{lemma}[The Hautus test]  \label{lem-Hautus}The following properties are equivalent for the pair $(A, B)$:
\begin{itemize}
\item[a)] $(A, B)$ is controllable.
\item[b)] $\mbox{ rank } [\lambda I - A, B] = n$ for all $\lambda \in \mC$.
\item[c)] $\mbox{ rank } [\lambda I - A, B] = n$ 
for all eigenvalues $\lambda$ of $A$. 
\end{itemize}
\end{lemma}

\begin{remark} \label{rem-Hautus} \rm Recall that, for an invertible $n \times n$ matrix $T$,  
\be \label{lem-Hautus-p1}
\mbox{rank } [A, B] =   \mbox{rank } [A T, B]
\ee
and 
\be \label{lem-Hautus-p2}
\mbox{rank } [A, B] =   \mbox{rank } [T A, T B]. 
\ee

In fact 
$$
\mbox{rank } C = \dim \mbox{Col} C, 
$$
$$
\mbox{Col} [A, B] = \{Ax + Bu; x \in \mR^n, u \in \mR^p \}
$$
$$
\mbox{Col} [A T, B] = \{ATx + Bu; x \in \mR^n, u \in \mR^p \}
$$
$$
\mbox{Col} [T A, TB] = \{TAx + TBu; x \in \mR^n, u \in \mR^p \} = T \mbox{Col} [A, B]. 
$$
So there is no difference in the Hautus test, if one changes the condition $\mbox{ rank } [\lambda I - A, B] = n$ by the condition $\mbox{ rank } [A - \lambda I, B] = n$.
\end{remark}

\begin{proof} We first make two observations. The first one is: if  $(A, B)$ is controllable then $\mbox{rank } (A, B) = n$ since otherwise, by a change of variables using the Gaussian elimination method, one has 
$$
\widetilde x ' = \widetilde A \widetilde x + \widetilde B \widetilde u 
$$
with the last line of $ (\widetilde A, \widetilde B)$ is $0$ and the system is not controllable. This fact can be also seen from \Cref{lem-Kalman-cond} and \Cref{rem-Hautus}. 

The second observation is that if $x$ is a solution of $x' = Ax + Bu$ then 
$$
(e^{-\lambda t} x)' = (A - \lambda I) (e^{-\lambda t} x) + B (e^{-\lambda t} u). 
$$
Hence if $(A, B)$ is controllable then $(A - \lambda I, B)$ is also controllable. 

Using these two observations, we are ready to give the proof of Hautus' lemma. 

\medskip 
It is clear that $a) \Rightarrow b)$ by the two observations. 

\medskip 

It is clear that $b) \Rightarrow c)$ and also $c) \Rightarrow b)$.  Thus $b)$ is equivalent to $c)$. 

\medskip 
It remains to prove that $b) \Rightarrow a)$ by contradiction as follows. Assume that $(A, B)$ is not controllable. Then by the Kalman decomposition \Cref{lem-Kalman-cond}, there exists an invertible $n \times n$ matrix $T$ such that 
$$
\widetilde A = T^{-1}  A T,  \quad \widetilde B =T^{-1}  B, 
$$
and 
\begin{equation}
\widetilde A  =  \left(\begin{array}{cc} A_1 & A_2 \\
0 & A_3
\end{array}\right),  \quad \widetilde B  =  \left(\begin{array}{c} B_1  \\
0 
\end{array}\right), 
\end{equation}
with $A_1 \in \mR^{r \times r}$ and $B_1 \times \mR^{r \times m}$ where $r = \mbox{rank } M$.   It follows that 
 \begin{equation}
[A - \lambda I, B] = [T (\widetilde A - \lambda I ) T^{-1}, T\widetilde B]
\end{equation}
This yields  
$$
\mbox{rank }[A - \lambda I, B] \mathop{=}^{ \eqref{lem-Hautus-p2}} \mbox{rank } [(\widetilde A - \lambda I) T^{-1},  \widetilde B] \mathop{=}^{ \eqref{lem-Hautus-p1}} \mbox{rank }[ \widetilde A - \lambda I,  \widetilde B]. 
$$
Note that 
\begin{equation}
  (\widetilde A - \lambda I,  \widetilde B) = \left(\begin{array}{ccc} A_1 - \lambda I & A_2  & B_1 \\
0 &A_3 -  \lambda I & 0
\end{array}\right).
\end{equation}
It follows that if  $\lambda$ is an eigenvalue of $A_3$, we have 
$$
\mbox{rank } [\widetilde A - \lambda,  \widetilde B] \le n-1: 
$$
which yields a contradiction. 
\end{proof}

\section{Observability}

Consider the linear control system 
\be\label{Ob-sys1}
x' = A x + Bu. 
\ee
An observation is given by the relation, for some $m \times n$ matrix $C$ ($m \ge 1$),  
\be\label{Ob-observation}
y =  C x. 
\ee

\begin{remark} \rm 
From practical viewpoint, $y$ is the variable one can observe not the full state $x$. 
\end{remark}

One next introduces the notation of observability. 

\begin{definition} System \eqref{Ob-sys1} and \eqref{Ob-observation}
is said to be observable if for $x_0 \neq z_0$, there exists $\tau > 0$ and  $u \in L^1((0, \tau); \mR^p)$ such that 
$$
C x(\tau) \neq C z(\tau), 
$$
where 
$$
x' = A x + Bu \mbox{ in } (0, \tau), \quad x(0) = x_0, \quad \mbox{ and } \quad  z' = A z + Bu \mbox{ in } (0, \tau), \quad z(0) = z_0
$$
If for a given $T > 0$ and for arbitrary $x \neq z$,  there exists $\tau \in  (0,T]$  with the above property, then System \eqref{Ob-sys1} and \eqref{Ob-observation} is said to be observable at time $T$. 
\end{definition}

We have the following property on the observability. 

\begin{proposition} \label{pro-Obs} System \eqref{Ob-sys1} and \eqref{Ob-observation} is observable if and only if for $x_0 \neq 0$, there exists $\tau > 0$ such that 
$$
C x(\tau) \neq 0, 
$$
where 
$$
x' = A x \mbox{ in } (0, \tau), \quad x(0) = x_0. 
$$
Consequently, system \eqref{Ob-sys1} and \eqref{Ob-observation} is observable in time $T>0$ if and only if 
for $x_0 \neq 0$, there exists $\tau \in (0, T]$ such that 
$$
C x(\tau) \neq 0, 
$$
where $x$ is the solution of the system 
$$
x' = A x \mbox{ in } (0, \tau), \quad x(0) = x_0. 
$$
\end{proposition}

It follows from \Cref{pro-Obs} that the observability depends only on $A$ and $C$. 

\medskip 
Set, for $T>0$,  
$$
R_T = \int_0^T e^{t A^\top} C^\top C e^{t A} \, dt.  
$$

We have the following result. 

\begin{theorem} \label{thm-Obs}
The following conditions are equivalent. 
\begin{itemize}
\item[i)] System \eqref{Ob-sys1} and \eqref{Ob-observation}   is observable.
\item[ii)] System \eqref{Ob-sys1} and \eqref{Ob-observation}   is observable at a given time $T > 0$. 
\item[iii)] $R_T$ is invertible for some $T> 0$. 
\item[iv)] $R_T$ is invertible for all $T> 0$. 
\item[v)] $\mbox{rank} \left[\begin{array}{c}
C\\
CA\\
\vdots\\
C A^{n-1}
\end{array} \right] = n$.  
\end{itemize}
\end{theorem}

\begin{proof} The proof is divided into several steps. 

\medskip  
\noindent {\it Step 1: $ii) \Rightarrow iii)$} We have 
$$
\langle R_T x, x \rangle = \int_0^T |C e^{t A} x|^2 \, dt.
$$
Assume that $R_T$ is not invertible. Then $\exists x \in \mR^n \setminus \{0 \}$ such that 
$$
R_T x = 0. 
$$
This impies 
$$
\int_0^T |C e^{t A} x|^2 \, dt = 0. 
$$
It follows that  
$$
C e^{t A} x = 0 \mbox{ for } t \in [0, T].  
$$
We have a contradiction. \footnote{One can check by the same arguments that $iii) \Rightarrow ii)$.}

\medskip 
\noindent {\it Step 2: $iii) \Rightarrow iv)$} Assume that $R_{T'}$ is not invertible for some $T' > 0$. There are some $x \in \mR^n \setminus \{0 \}$ such that  
$$
C e^{t A} x = 0 \mbox{ in } [0, T']. 
$$
This implies, by the analyticity of the function $C e^{\cdot A} x$, 
$$
C e^{t A} x = 0 \mbox{ for } t \in \mR. 
$$
We derive that 
$$
R_{T}x = \int_0^T e^{t A^\top }C^\top C e^{t A} x \, dt = 0 \mbox{ for all } T > 0. 
$$
We have a contradiction. \footnote{It is clear that $iv) \Rightarrow iii)$.}

{\it Step 3: $iv) \Rightarrow v)$} Assume that the rank is less than $n$. There exists $x \in \mR^n \setminus \{0 \}$ such that 
$$
C A^k x = 0 \mbox{ for } 0 \le k \le n-1. 
$$
By the Caley-Hamilton theorem, 
$$
C A^k x = 0 \mbox{ for } k \ge 0.  
$$
This implies that 
$$
f^{(k)}(0) = 0 \mbox{ for } k \ge 0 \mbox{ where }
f(t) = C e^{t A}. 
$$
By the analyticity of $f$, we have 
$$
f(t) = 0 \mbox{ in } \mR. 
$$
This implies that 
$$
R_T(x) = \int_0^T e^{t A^\top }C^\top C e^{t A} x \, dt = 0. 
$$

{\it Step 4: $v) \Rightarrow ii)$} It is clear that $R_T$ is invertible for all $T>0$. This implies, for $x \in \mR^n \setminus \{0\}$, 
$$
C e^{t A}x \neq 0 \mbox{ for some } t> 0.  
$$

The proof is complete. 
\end{proof}

\begin{remark} \rm In what follows, we often say the pair $(A, C)$ is observable instead of saying that System \eqref{Ob-sys1} and \eqref{Ob-observation}  is observable. This completely makes sense by \Cref{pro-Obs} or 
\Cref{thm-Obs}. 
\end{remark}

As a consequence of \Cref{thm-LA-Kalman} and \Cref{thm-Obs}, one obtains the following duality between the observability and the controllability. 

\begin{corollary} $(A, C)$ is observable if and only if  $(A^\top, C^\top)$ is controllable. 
\end{corollary}

\begin{remark} \rm The materials of this section are standard, see, e.g., \cite{Sontag98,Zabczyk08}. 
\end{remark}

\section{Stability} \label{sect-stability}

In this section, we discuss the stability of the system $x' = Ax$ and explore some roles of the concept of observability.  Consider the linear system 
\be \label{SS-sys1}
x'=Ax \mbox{ and } x(0) = x_0 \in \mR^n. 
\ee
The solution is then given by 
$$
x(t) = e^{tA} x_0 \mbox{ for } t \ge 0. 
$$
System \eqref{SS-sys1} is called stable if 
$$
e^{tA} x_0 \to 0 \mbox{ as } t \to + \infty \mbox{ for all } x_0 \in \mR^n. 
$$
Instead of saying that \eqref{SS-sys1} is stable one often says that the matrix $A$ is stable. 


\medskip 
We recall the following well-known result, whose proof can be found in any textbook of differential equations. 
 
\begin{theorem} \label{thm-ODE}The following conditions are equivalent: 
\begin{itemize}
\item[i)] $A$ is stable. 

\item[ii)] $e^{tA} x_0$ exponentially decays as $t \to + \infty$ for all $x_0 \in \mR^n$. 

\item[iii)] $\omega(A) = \sup \Big\{ \Re \lambda, \lambda \mbox{ is an eigenvalue of $A$} \Big\} < 0$ \footnote{Here $\Re \lambda$ denotes the real part of $\lambda$.}.  

\item[iv)] $\dsp \int_0^\infty |e^{tA} x_0|^2 \, dt < + \infty \mbox{ for all $x_0 \in \mR^n$}.$ 
\end{itemize}
\end{theorem}

An important role in the stability of the system $x' = Ax$ is played by the so-called matrix Lyapunov equation of the form
\be\label{SS-Lyapunov-equation}
A^\top Q + Q A = - R
\ee
where $R$ is an $n \times n$ non-negative, symmetric matrix, and $Q$ is an (unknown) $n \times n$ symmetric matrix. In this spirit, the following result is given.

\begin{lemma} \label{lem-ODE} Consider the system $x' = Ax$. Let $R$ be symmetric and positive, and assume that there exists a solution $Q$ symmetric and positive of \eqref{SS-Lyapunov-equation}. Then $ii)$ of \Cref{thm-ODE} holds. 
\end{lemma}  

\begin{proof} Given $x_0 \in \mR^n$, let $x$ be the unique solution of the system 
$$
x' = A x \mbox{ in } [0, + \infty) \quad \mbox{ and } \quad x(0) = x_0. 
$$
Set 
$$
V(t) = \langle Q x(t), x(t) \rangle \mbox{ for } t \ge 0. 
$$
We have 
\begin{multline*}
\frac{d}{dt} V(t) =  \langle Q x'(t), x (t) \rangle + \langle Q x(t), x' (t) \rangle = \langle Q A x (t), x (t) \rangle + \langle Q x(t), A x(t) \rangle \\[6pt]
= \langle Q A x (t), x (t) \rangle  + \langle A^\top Q x(t), x(t) \rangle, 
\end{multline*}
which yields, by \eqref{SS-Lyapunov-equation},  
\be \label{lem-ODE-p1}
\frac{d}{dt} V(t) = - \langle R x(t), x(t) \rangle. 
\ee

Since $R$ is symmetric and positive, there exists $c_1 > 0$ such that 
$$
\langle R z, z \rangle \ge c_1 \langle Q z, z \rangle \mbox{ for all } z \in \mR^n. 
$$
It follows that 
$$
\frac{d}{dt} V(t) \le - c_1 V(t) \mbox{ for } t \ge 0. 
$$
Thus 
$$
V(t) \le V(0) e^{-c_1 t} \mbox{ for } t \ge 0. 
$$
Since $Q$ is symmetric and positive, there exists $c_2 > 0$ such that 
$$
\langle Q z, z \rangle \ge c_2 \langle z, z \rangle \mbox{ for all } z \in \mR^n. 
$$
We derive that 
$$
\| x(t)\|^2 \le c_2^{-1} V(0) e^{-c_1 t}. 
$$
Thus $A$ is stable. The proof is complete. 
\end{proof}

Here is a more subtle version of \Cref{lem-ODE} which involves the observability notion.

\begin{theorem} \label{SS-thm1} We have
\begin{itemize}
\item[i)] Assume that the pair $(A, C)$ is observable and denote $R = C^\top C$. If there exists a nonnegative definite matrix $Q$ satisfying \eqref{SS-Lyapunov-equation} then the matrix $A$ is stable \footnote{As seen in the proof,  $Q$ is positive by using the fact that $(A, C)$ is observable.}.

\item[ii)] If the matrix $A$ is stable in the sense of $iv)$ of \Cref{thm-ODE}, then given an arbitrary symmetric matrix $R$, equation \eqref{SS-Lyapunov-equation} has exactly one symmetric solution $Q$. This solution is positive (nonnegative) definite if the matrix $R$ is positive (nonnegative) definite.
\end{itemize}
\end{theorem}

\begin{proof} The proof explore the ideas given in \Cref{lem-ODE} and the property of observability. The proof is divided into two steps. 

\medskip 
\noindent {\it Step 1}: Proof of $i)$. Denote $S(t) = e^{t A}$ for $t \in \mR$.  As in the proof of \eqref{lem-ODE-p1}, we have, by \eqref{SS-Lyapunov-equation},  
\be \label{SS-thm1-p1-**}
\frac{d}{dt} S(t)^\top Q S(t) = - S(t)^\top  R S(t).  
\ee
Set 
\be \label{SS-thm1-p0}
R_T = \int_0^T S(t)^\top  R S(t) \, dt \mathop{=}^{R = C^\top C} \int_0^T S(t)^\top  C^\top C S(t) \, dt. 
\ee
Integrating \eqref{SS-thm1-p1-**} from $0$ to $T$ and taking into account the fact $S(0) = I$, one obtains  
\be \label{SS-thm1-p1}
Q  = R_T + S(T)^\top Q S(T).   
\ee
Since $(A, C)$ is observable, it follows from \Cref{thm-Obs} that 
\be \label{SS-thm1-p2}
R_T > 0. 
\ee
Combining \eqref{SS-thm1-p1} and \eqref{SS-thm1-p2} and using the fact that $Q \ge 0$, we derive that 
\be \label{SS-thm1-p3}
Q > 0. 
\ee

Given $x_0 \in \mR^n$, let $x$ be the solution of the system 
$$
x' = A x \mbox{ in } \mR \quad \mbox{ and }  \quad x(0) = x_0. 
$$
We have, as in the proof of \eqref{lem-ODE-p1},  
\be \label{SS-thm1-p4}
\frac{d}{dt} \langle Q x(t), x(t) \rangle  = - \langle R x(t), x(t) \rangle = - \|C x(t) \|^2. 
\ee
Since $Q$ is positive, we derive that $\Re \lambda \le 0$ for all eigenvalue $\lambda$ of $A$. We next show by contradiction that 
$$
\Re \lambda < 0 \mbox{ for all eigenvalues $\lambda$ of $A$.}
$$
Assume that there exists an eigenvalue $\lambda$ of $A$ such that $\lambda = i \omega$ for some $\omega \in \mR$. We first note that $\omega \neq 0$ since otherwise $\exists \xi \in \mR^n$ such that $A \xi = 0$ and thus $S(t) \xi = \xi$ for $t \in \mR$ and $C \xi = 0$. This contradicts the observability property of $(A, C)$, see \Cref{pro-Obs}.

We thus have $\omega \neq 0$. Let $z_0 \in \mC^n \setminus 0$ be an corresponding eigenvalue. We have 
$$
z_0 = x_0 + i y_0 \mbox{ where } x_0, y_0 \in \mR^n. 
$$
Since
$$
A(x_0 + i y_0) = i \omega (x_0 + i y_0), 
$$
it follows that 
$$
A x_0 = - \omega y_0 \quad \mbox{ and } \quad A y_0 = \omega x_0. 
$$
Set  
$$
\xi(t) = \cos (\omega t) x_0 - \sin (\omega t) y_0. 
$$
Then 
$$
\xi'(t) = A \xi(t). 
$$
The periodicity of $\xi$ and \eqref{SS-thm1-p4} imply that 
$$
\langle Q \xi(t), \xi(t) \rangle \mbox{ is constant for $t \ge 0$}. 
$$
This implies, again by \eqref{SS-thm1-p4},   
$$
C \xi(t) = 0 \mbox{ for } t \ge 0. 
$$
We have a contradiction with the observability of $(A, C)$, see \Cref{pro-Obs}. 

\medskip 
\noindent{\it Step 2}: Proof of $ii)$. The analysis is motivated from \eqref{SS-thm1-p1}. Let $R$ be symmetric. Set 
\be \label{SS-thm1-def-Q}
Q   = \int_0^\infty S(t)^\top R S(t) \, dt.  
\ee
Since $A$ is stable in the sense of $iv)$ of \Cref{thm-ODE}, it follows that $Q$ is well-defined. It is clear that $Q$ is symmetric since $R$ is symmetric. We have
\begin{multline} \label{SS-thm1-p00}
A^\top Q + Q A =  \int_0^\infty \Big( A^\top S(t)^\top R S(t) + S(t)^\top R S(t) A  \Big) \, dt \\[6pt]
= \int_0^\infty \left( \frac{d}{dt}S(t)^\top R S(t) + S(t)^\top R \frac{d}{dt} S(t)  \right) \, dt =  \int_0^\infty \frac{d}{dt} S(t)^\top R S(t) \, dt = -R. 
\end{multline}
which implies the existence of a solution of \eqref{SS-Lyapunov-equation}. To prove the uniqueness, we note that if 
$$
A^\top Q + Q A =0, 
$$
then
$$
\frac{d}{dt} S(t)^\top Q S(t) =0. 
$$
Since 
$$
\lim_{t \to + \infty}S(t)^\top Q S(t) =0, 
$$
it follows that 
$$
S(0)^\top Q S(0) =0. 
$$
Since $S(0) = I$, we derive that $Q = 0$. The uniqueness follows. 

It is clear to see that $Q$ given by \eqref{SS-thm1-def-Q} is non-negative (resp. positive) if $R$ is non-negative (resp. positive). 

The proof is complete. 
\end{proof}

\begin{remark} \rm \Cref{lem-ODE} and  part ii) of \Cref{SS-thm1} have an interesting counter part in the infinite dimensional setting, see, e.g., \cite{BDDM07}. 

\end{remark}

Here is an interesting consequence of \Cref{SS-thm1}. 

\begin{corollary}\label{SS-cor1} Consider the system  \eqref{Ob-sys1} and \eqref{Ob-observation}.  Assume that $(A, C)$ is observable and 
\be \label{SS-cor1-ass1}
C e^{t A} x_0 \to 0 \mbox{ as } t \to + \infty \mbox{ for all } x_0 \in \mR^n. 
\ee
Then $A$ is stable, and consequently 
\be
e^{t A} x_0 \to 0 \mbox{ exponentially as } t \to + \infty \mbox{ for all } x_0 \in \mR^n. 
\ee
\end{corollary}

\begin{proof} Let $\beta > 0$ (arbitrary). Set $S_\beta(t) = e^{ t(A - \beta I)}$ for $t \in \mR$ \footnote{Recall that, for $A_1$ and $A_2$ are two $n\times n$ matrices such that $A_1 A_2 = A_2 A_1$, then $e^{A_1 + A_2} = e^{A_1} e^{A_2}$.}, and 
$$
Q_\beta = \int_0^\infty S_\beta(t)^\top C^\top C S_\beta(t) \, dt.  
$$
The matrix $Q_\beta$ is well-defined by \eqref{SS-cor1-ass1}. We have, similar to \eqref{SS-thm1-p00},  
\be \label{SS-cor1-p1}
(A - \beta I)^\top Q_\beta + Q_\beta (A - \beta I) = - C^\top C. 
\ee
It follows from $v)$ of \Cref{thm-Obs}, for $\beta$ sufficiently small, that $(A - \beta I, C)$ is observable. Applying  $i)$ of \Cref{SS-thm1}, we derive from \eqref{SS-cor1-p1} that, for $\beta > 0$ small,  
$$
A - \beta I \mbox{ is stable}. 
$$
This implies, for $\beta > 0$ sufficiently small. 
$$
\Re \lambda - \beta < 0
$$
for all eigenvalue $\lambda$ of $A$. 
This implies 
$$
\Re \lambda  \le 0. 
$$

We next prove that 
$$
\Re \lambda < 0, 
$$
by contradiction. The proof is similar to the one of the last part of Step 1 of the proof of \Cref{SS-thm1}.  Assume that $\lambda = i \omega$ for some $\omega \in \mR$. Since $(A, C)$ is observable, $\omega \neq 0$. Let $x_0, y_0 \in \mR^n$ be such that $(x_0, y_0) \neq 0$, 
$$
A (x_0 + i y_0) = i \omega (x_0 + i y_0)
$$
it follows that 
$$
A x_0 = - \omega y_0 \quad \mbox{ and } \quad A y_0 = \omega x_0. 
$$
Set  
$$
\xi(t) = \cos (\omega t) x_0 - \sin (\omega t) y_0. 
$$
Then 
$$
\xi'(t) = A \xi(t). 
$$
It follows that 
$$
C \xi(t) \to 0. 
$$
Since $\xi$ is periodic, one derives that 
$$
C \xi = 0 \mbox{ for }  t \ge 0. 
$$
This contradicts the observability of $(A, C)$. 
\end{proof}

\begin{remark} \rm The materials presented in this section are standard, see, e.g.,  \cite{Sontag98}. The presentation given here follows from \cite{Zabczyk08}. 
\end{remark}

\section{Detectability and dynamic observers}

Here is the notion of detectability. 

\begin{definition} A pair of matrices $(A, C) \in \mR^{n \times n} \times \mR^{m \times n}$ is said to be detectable if there exists an $n \times m$ matrix $L$ such that the matrix $A + LC$ is stable. 
\end{definition}

The following result gives a relation between observability and detectability. 

\begin{proposition} Assume that $(A, C)$ is observable. Then $(A, C)$ is detectable. 
\end{proposition}

\begin{proof} Since $(A, C)$ is observable, $(A^\top, C^\top)$ is controllable. There thus exists $F$ such that $A^\top + C^\top F$ is stable by Pole-Shifting Theorem \Cref{thm-PS} proved later. Equivalently, $A + F C$ is stable. One chooses then $L = F$. 
\end{proof}

We illustrate the importance of the concept of detectability by discussing the dynamical observer introduced by Luenberger \cite{luenberger1971introduction}. To this end, one needs the concept of stabilizability. 

\begin{definition} A pair $(A, B) \in \mR^{n \times n} \times \mR^{p \times n}$ is said to be stabilizable if there exists an $n \times p$ matrix $K$ such that the matrix $A + BK$ is stable. 
\end{definition}

The following is the main result of this section. 

\begin{theorem} \label{thm-DO} Consider the control system 
$$
x' = A x + Bu   
$$
with the observation 
$$
y = Cx.
$$
Assume that the pair $(A, C)$ is detectable and let $L$ be a matrix such that the matrix $A + LC$ is stable. Then the system 
$$
\hx' = (A + LC) \hx - L y + Bu \\[6pt]
$$
defines a dynamical observer, i.e., $ \hx - x \to 0$ exponentially as $t \to + \infty$. Assume in addition that the pair $(A, B)$ is stabilizable and let $K$ be a matrix such
that $A + BK$ is stable. Then, for the control $u = K \hx$, it holds 
$$
x(t) \to 0 \mbox{ as } t \to + \infty \mbox{ exponentially}. 
$$
\end{theorem}

\begin{remark} \rm It follows from \Cref{thm-DO} that the system $x' = Ax + Bu$ with the observation $y = Cx$ can be exponentially stabilized for which the feedback uses only the information of the observer $y$ and the variables $\hat x$, which depends only on $y$ and $u$. 
\end{remark}

\begin{proof} We have
$$
(\hx - x)' = (A + LC) \hx - L C x + Bu - (Ax + Bu) = (A + LC) (\hx - x). 
$$
Since $A + LC$ is stable, we have 
$$
\|(\hx - x)(t) \| \le c_1 e^{-\lambda_1 t} \|(\hx - x)(0) \| \mbox{ for } t \ge 0, 
$$
for some positive constants $c_1$ and $\lambda_1$. The first assertion follows. 

We next deal with the second assertion. We have 
$$
\hx' = (A + L C) \hx - L C x + B K \hat x = (A + BK) \hx + LC (\hx - x). 
$$
It follows that 
\be \label{thm-DO-p1}
\hx(t) = e^{(A + BK) t} \hx(0) + \int_0^t  e^{(A + BK) (t-s)} L C (\hx - x)(s) \, ds. 
\ee
Since $A + BK$ is stable, there exist some positive constants $c_2$ and $\lambda_2$ such that, for $ \xi \in \mR^n$,  
$$
\| e^{(A + BK) t} \xi \| \le c_2 e^{-\lambda_2 t} \| \xi\| \mbox{ for } t \ge 0. 
$$
It follows from the first assertion that 
\begin{multline} \label{thm-DO-p2}
\|e^{(A + BK) t} \hx(0)\| + \int_0^t  \left|  e^{(A + BK) (t-s)} L C (\hx - x)(s) \right|\, ds \\[6pt]
\le c_2 e^{-\lambda_2 t} \| \hx(0)\| + c_3 \| \hx(0) - x(0)\|  \int_0^t e^{-\lambda_2(t-s)} e^{- \lambda_1 s} \, ds,  
\end{multline}
for some positive constant $c_3$. Since 
$$
 \int_0^t e^{-\lambda_2(t-s)} e^{- \lambda_1 s}  \le c_4 e^{-\lambda t} 
$$
for some positive constants $c_4$ and $\lambda$, the second assertion follows from \eqref{thm-DO-p1} and \eqref{thm-DO-p2}. 
\end{proof}

\section{Pole shifting theorem}

We begin this section with the following definition. 

\begin{definition} \label{def-equivalenceAB} Let $(A,B)$ and $(\tA, \tB)$ be two pairs of matrices in $\mR^{n \times n} \times \mR^{n \times p}$.  Then (A,B) is
similar to $(\tA, \tB)$ and this is denoted by $(A, B) \sim (\tA, \tB)$ if there exists $T \in \mR^{n \times n}$ invertible such that 
\be \label{def-equivalenceAB-tAtB}
T^{-1} A T = \tA \quad \mbox{ and } \quad T^{-1} B = \tB. 
\ee
\end{definition}

\begin{remark} \rm The motivation for this definition is as follows. Consider the control system 
$$
x' = A x + B u. 
$$
Let $T \in \mR^{n \times n}$ be invertible and set $\tx = T^{-1} x$. We then have 
$$
(T \tx)' = A T \tx + B u. 
$$
This is equivalent to the fact 
$$
\tx' = T^{-1} A T \tx + T^{-1} B u = \tA \tx + \tB u,  
$$
where $(\tA, \tB)$ is defined by \eqref{def-equivalenceAB-tAtB}. Thus the control system given by $(\tA, \tB)$ is obtained from the control system given by $(A, B)$ after a change of variables. 
\end{remark}

Concerning \Cref{def-equivalenceAB}, we have the following properties. 

\begin{lemma} The following facts hold for the relation $\sim$:
\begin{itemize}
\item  It is an equivalence relation on $\mR^{n \times n} \times \mR^{n \times p}$.

\item  Assume that  $(A, B) \sim (\tA, \tB)$. Then $(A, B)$ is controllable if and only if $(\tA, \tB)$ is. 
\end{itemize}
\end{lemma}

In what follows in this section, when $p=1$, we often denote $B$ by $b$. 

\begin{lemma} \label{lem-AB-hAhB} Let $(A, b)  \in \mR^{n \times n} \times \mR^{n \times 1}$ and assume that the characteristic polynomial $\chi_A$ of $A$ is given by 
$$
\chi_A(s) = \det (sI - A) = s^n - \alpha_n s^{n-1} - \cdots - \alpha_2 s - \alpha_1. 
$$
Define 
$$
\hA = \left(\begin{array}{ccccc}
0 & 0 & \dots & 0 & \alpha_1 \\
1 & 0 & \cdots & 0 & \alpha_2 \\ 
0 & 1 & \cdots & 0 & \alpha_3 \\
\cdots & \cdots & \cdots & \cdots & \cdots \\
0 & 0 & \cdots &  1 & \alpha_n \\
\end{array} \right) \quad \mbox{ and } \quad  \hat b = \left(\begin{array}{c}
1\\
0 \\
0 \\
\cdots \\
0 
\end{array}\right).
$$
Then 
$$
A \cR(A, b) = \cR(A, b) \hat A \quad \mbox{ and } \quad b = \cR (A, b) \hat b. 
$$
where $\cR (A, b) = (b, Ab, \dots, A^{n-1}b)$. 
Consequently, $(A, b)$ is controllable if and only if it is similar to $(\hat A, \hat b)$. 
\end{lemma}

\begin{proof}  The first part is an immediate consequence of the Cayley-Hamilton Theorem. The second part is clear since $(\hA, \hb)$ is controllable by the Kalman rank condition.
\end{proof}

\begin{remark} \rm It is worth noting that 
\be \label{lem-AB-hAhB-det}
\det (sI - \hat A) = \det (s I - A). 
\ee
This can be proved by recurrence as follows. The identity is clear for $n =1$. Assume $n \ge 2$ and the identity holds for $n-1$. We have 
\begin{multline*}
\det \left(\begin{array}{ccccc}
s & 0 & \dots & 0 & -\alpha_1 \\
-1 & s & \cdots & 0 & -\alpha_2 \\ 
0 & -1 & \cdots & 0 &  - \alpha_3 \\
\cdots & \cdots & \cdots & \cdots & \cdots \\
0 & 0 & \cdots &  -1 & s - \alpha_n \\
\end{array} \right) \\[6pt]= s \det \left(\begin{array}{ccccc}
 s & \cdots & 0 & -\alpha_2 \\ 
 -1 & \cdots & 0 &  - \alpha_3 \\
 \cdots & \cdots & \cdots & \cdots \\
 0 & \cdots &  -1 & s - \alpha_n \\
\end{array} \right) + (-1)^{n+2} \alpha_1  \det  \left(\begin{array}{ccccc}
-1 & s & \cdots & 0  \\ 
0 & -1 & \cdots & 0  \\
\cdots & \cdots & \cdots & \cdots \\
0 & 0 & \cdots &  -1  \\
\end{array} \right). 
\end{multline*}
We have, by the recurrence,  
$$
s \det \left(\begin{array}{ccccc}
 s & \cdots & 0 & -\alpha_2 \\ 
 -1 & \cdots & 0 &  - \alpha_3 \\
 \cdots & \cdots & \cdots & \cdots \\
 0 & \cdots &  -1 & s - \alpha_n \\
\end{array} \right)  = s (s^{n-1} - \alpha_n s^{n-2} - \dots - \alpha_2). 
$$ 
Since 
$$
\det  \left(\begin{array}{ccccc}
-1 & s & \cdots & 0  \\ 
0 & -1 & \cdots & 0  \\
\cdots & \cdots & \cdots & \cdots \\
0 & 0 & \cdots &  -1  \\
\end{array} \right)  = (-1)^{n-1},  
$$ 
assertion \eqref{lem-AB-hAhB-det} follows.  
\end{remark}

The system $(\hA, \hb)$ is sometimes called the controllability form of $(A, b)$. From the above lemma and the fact that $(\hA, \hb)$ depends only on $\chi_A$ and since $\sim$ is an equivalence relation, we reach the following.

\begin{corollary} \label{cor-similarity} Let $(A, b)$ and $(\tA, \tb)$ be in $\mR^{n \times n} \times \mR^{n \times 1}$. 
Assume that the pairs $(A, b)$ and $(\tA, \tb)$ are controllable. Then $(A, b)$ is similar to $(\tA, \tb)$ if and only if $\chi_A = \chi_{\tA}$. 
\end{corollary}

\begin{remark} \rm The condition on the controllability of  $(A, b)$ and $(\tA, \tb)$ is necessary. 
\end{remark}

\begin{proof} 
We have 
$$
(A, b) \sim (\hat A, e_1) \quad \mbox{ and } \quad (\tA, b) \sim (\hat \tA, e_1).
$$ 
Thus if  the characteristic polynomials of $A$ and $\tA$ are the same then $ (\hat A, e_1)  =  (\hat \tA, e_1) $. Hence $(A, b) \sim (\tA, \tb)$. 

Assume now that $(A, b) \sim (\tA, \tb)$. Then there exists an $(n \times n)$  invertible matrix $T$ such that 
\be \label{cor-PS-p1}
\hat A = T^{-1} \tA T. 
\ee
Therefore, the characteristic polynomials of $A$ and $\tA$ are the same. 
\end{proof}

\begin{definition} The controller form  associated to the pair  $(A, b)$ is the pair 
$$
A^{\sharp} =  \hA^\top =  \left(\begin{array}{ccccc}
0 & 1 & \dots & 0 & 0 \\
0 & 0 & 1 & 0 & 0 \\ 
\cdots & \cdots & \cdots & \cdots & \cdots \\
0 & 0 & 0 & \cdots &  1 \\
\alpha_1 & \alpha_2 & \cdots &  \alpha_{n-1} & \alpha_n \\
\end{array} \right) \quad \mbox{ and } \quad  b^\sharp = \left(\begin{array}{c}
0\\
0 \\
0 \\
\cdots \\
1 
\end{array}\right),
$$
where $s^n - \alpha_n s^{n-1} - \cdots - \alpha_2 s - \alpha_1$ is the characteristic polynomial of $A$. 
\end{definition}

Concerning $(A^{\sharp}, b^{\sharp})$,  we have the following result. 
 
\begin{lemma} \label{lem-control-form} The pair $(A^\sharp, b^\sharp)$ is controllable and the characteristic polynomial of $A^\sharp$ is $s^n - \alpha_n s^{n-1} - \cdots - \alpha_2 s - \alpha_1$. 
\end{lemma}

\begin{proof} It is a direct consequence of the Kalman rang criterion applied to $(A^\sharp, b^\sharp)$. 
\end{proof}

\begin{proposition} The  pair $(A, b)$ is controllable if and only if it is similar to its controller form.
\end{proposition}

\begin{proof} Assume first that $(A, b)$ is controllable. Since $(A^\sharp, b^\sharp)$ is controllable by \Cref{lem-control-form}, it follows from \Cref{cor-similarity} that  $(A, b) \sim (\hat A, \hat b)$ and $(\hat A, \hat b) \sim (A^\sharp, b^\sharp)$. This implies $(A, b) \sim (A^\sharp, b^\sharp)$

Assume next that $(A, b) \sim (A^\sharp, b^\sharp)$. Since $(A^\sharp, b^\sharp)$ is controllable, we derive that $(A, b)$ is controllable. 
\end{proof}

Thus, up to a change of variables, single-input controllable linear, time-invariant systems are precisely those obtained from $n^{th}$ order constant coefficient scalar differential equations
$$
x^{(n)} - \alpha_n x^{(n-1)} - \dots - \alpha_1 x = u, 
$$
after a change of variables.  We shall be interested in studying the effect of constant linear feedback. The following lemma will be helpful in that regard.

\begin{lemma} For any pair $(A, B)$ and any $F \in \mR^{p \times n}$
$$
\cR(A + BF, B) = \cR(A, B).
$$
In particular, $(A + BF, B)$ is controllable if and only if $(A, B)$ is. 
\end{lemma}

\begin{proof} One way is to note that 
$$
\cR(A, B) : = \mbox{span } \Big\{ A^i b_j;  \; i \ge 0, 1 \le j \le m \Big\}. 
$$

The other way is to realize that if 
$$
x' = A x + Bu 
$$
and $u = Fx + v$ then 
$$
x' = (A + BF) x + Bv. 
$$
The controllablility of the two systems are equivalent. 
\end{proof}

The following definition is useful. 

\begin{definition} The $n^{th}$ degree monic polynomial is assignable for the
pair $(A, B)$ if there exists a matrix $F$ such that $\chi_{A + B F} = \chi$. 
\end{definition}

\begin{lemma} If $(A, B) \sim (\tA, \tB)$ then they can be assigned the same polynomials.
\end{lemma}

\begin{proof} Let $T$ be invertible such that $A = T^{-1} \tA T $ and $B = T^{-1} \tB$.
Given $F$, it holds that $A+BF = T^{-1} \tA T + T^{-1} \tB  F = T^{-1} \tA T + T^{-1} \tB  \tF T = T^{-1} (\tA + \tB \tF)$ with $\tF = FT^{-1}$. The conclusion follows. 
\end{proof}

Here is the main result of this section. 

\begin{theorem}[Pole-Shifting Theorem] \label{thm-PS} The pair $(A,B)$ is controllable if and only if every $n^{th}$ degree
monic polynomial can be assigned to it.
\end{theorem}

\begin{proof} The fact that  every $n^{th}$ degree
monic polynomial can be assigned to $(A, B)$ implies the controllability of the pair $(A, B)$ is a consequence of \Cref{lem-Kalman-cond} since one cannot make the action on the last $n-r$ components for the system $(\tA, \tB)$ where $r = \rank M$ otherwise. 

It remains to prove that if $(A, B)$ is controllable and it can assign with every monic polynomial.

The case $p=1$ follows from the controller form.  Now let $p \in \mN$ be arbitrary. Pick $v \in \mR^p$ such that $B v \neq 0$ and define $b = B v$. 

We next show that there exists an $F_1 \in \mR^{p \times n}$ s.t. $(A + B F_1, b)$ is controllable. Because the result has been established for the case when $p=1$, there exists  $f \in \mR^n$ such that the characteristic polynomial of $A + B(F_1 + v f) =  (A + B F_1) + b f $ is as desired,  and the conclusion follows by taking $ F = F_1 + v f$.

To establish the existence of $F_1$, we argue as follows. Let $\{x_1, \dots, x_k\}$ be any sequence of linearly independent elements in $\mR^n$ with $k$ as large as possible,  such that (with $x_0 :  = 0$ and $x_1 = Bv$)
$$
x_i - A x_{i-1} \in Col B. 
$$
We claim that $k = n$. Indeed, consider the space $V$ spanned by $\{x_1, \cdots, x_k \}$. By the maximality of $k$, 
$$
A x_k + B u \in V \mbox{ for } u \in \mR^p.  
$$
This implies, in particular, that 
$$
A x_k \in V. 
$$
It follows that 
$$
Col B \subset V. 
$$
Hence $A x_l \in V$ for $l = 1, \cdots, k-1$. 
So $V$ is an $A$-invariant subspace of $\mR^n$ containing $Col B$. It follows that $\cR(A, B) \subset V$. By the controllability, one has $k = n$. 

For $i = 1, \dots, n-1 = k-1$, let 
$$
F_1 x_i  = u_i, 
$$
where $u_i \in \mR^p$ is s.t.,  
$$
x_{i+1} - A x_{i}  = B u_{i}, 
$$
and we define $F_1 x_k$ arbitrary. Since 
$$
\cR(A + B F_1, x_1) = (x_1, \cdots, x_n), 
$$
the existence of $F_1$ with the desired properties is proved. 
\end{proof}

\begin{remark} \rm The materials presented in this section are standard, see, e.g., \cite{Sontag98,trelat2005controle}. 
\end{remark}

\begin{remark} \rm \Cref{thm-PS} is a powerful result on the stabilization in the finite-dimensional setting. Nevertheless, its generalization to the infinite dimensional setting is quite difficult and therefore, several other approaches are proposed to deal with the stabilization, see \Cref{sect-optimal} and \Cref{sect-Gramian}, and the references therein.  
\end{remark}

\section{The linear test on controllability} 

Consider the control system 
\begin{equation}\label{Sys-NL}
x' = f(x, u)
\end{equation}
where $x$ is a state in $\mR^n$ and $u$ is a control in $\mR^p$.  The following assumptions are used for $f$: 
\begin{equation}
\mbox{$f$ is of class $C^1$ in some open subset of $\mR^n \times \mR^p$}. 
\end{equation}

In this section, we study the local controllability along a trajectory. Here is the definition. 

\begin{definition} Let $t_0 < t_1$ and let $(\bar x, \bar u)$ be a $C^1$ trajectory of \eqref{Sys-NL} on $[t_0, t_1]$, i.e., $(\bar x, \bar u)$ is of class $C^1$ and $\bar x' = f(\bar x, \bar u)$ in $[t_0, t_1]$. The system \eqref{Sys-NL} is called locally controllable along $(\bar x, \bar u)$ if for all $\eps > 0$, there exists $\delta > 0$ such that for all $x_0 \in B(\bar x(t_0),  \delta)$ and $x_1 \in B(\bar x(t_1),  \delta)$, there exists  $C^1$ trajectory $(x, u)$ of \eqref{Sys-NL} on $[t_0, t_1]$
verifying 
$$
\|u - \bar u \|_{C^1([t_0, t_1])} \le \eps, \quad x(\tau) = x_0, \quad x(\sigma) = x_1. 
$$
\end{definition}

An important part in the understanding of the local controllability property is to understand the linearized system. Consider the system \eqref{Sys-NL} of class $C^1$. Let $(\bar x, \bar u)$ be a trajectory on $[t_0, t_1]$ and consider $(\bar x + \delta x, \bar u + \delta u)$ a trajectory nearby on $[t_0, t_1]$. We formally have  
$$
(\bar x + \delta x)' = f (\bar x + \delta x, \bar u + \delta u)  = f_x (\bar x, \bar u) \delta x + f_u(\bar x, \bar u) \delta u + o(\delta x, \delta u) \mbox{ in } [t_0, t_1].  
$$
Thus the linearized system is given by 
\begin{equation}\label{Sys-NL-L}
x'(t) = A(t) x(t) + B(t) u(t) \mbox{ in } [t_0, t_1], 
\end{equation}
where
\begin{equation}\label{Sys-NL-LAB}
A(t) = f_x(\bar x(t), \bar u(t)) \mbox{ and } B(t) =  f_u(\bar x(t), \bar u(t)) \mbox{ in } [t_0, t_1]. 
\end{equation}
Since $f$ is of class $C^1$, it follows that $A$ and $B$ are continuous. 

\begin{example} Consider the pendulum problem. We have 
$$
f(x, u) = \left( \begin{array}{c}
x_2 \\
- \sin x_1 + u
\end{array}\right).  
$$
The linearized system is 
$$
x' = A(t)  x + B(t) u, 
$$
where 
$$
A(t) =  \left( \begin{array}{cc}
0 & 1  \\
- \cos x_1(t) & 0
\end{array}\right) \quad \mbox{ and } \quad 
B(t) =  \left( \begin{array}{c}
0 \\
1
\end{array}\right). 
$$
If one considers the equilibrium $(x, u) = (0, 0)$, then 
$$
A(t) =  \left( \begin{array}{cc}
0 & 1  \\
-1 & 0
\end{array}\right) \quad \mbox{ and } \quad 
B(t) =  \left( \begin{array}{c}
0 \\
1
\end{array}\right). 
$$
if one considers the equilibrium $(x, u) = (\pi, 0)$, then 
$$
A(t) =  \left( \begin{array}{cc}
0 & 1  \\
1 & 0
\end{array}\right) \quad \mbox{ and } \quad 
B(t) =  \left( \begin{array}{c}
0 \\
1
\end{array}\right). 
$$
\end{example}

Here is the main result of this section. 

\begin{theorem} \label{thm-LinearTest} Let $t_0 < t_1$. Assume that $f$ is of class $C^1$ and let $(\bar x, \bar u)$ be a $C^1$-trajectory of \eqref{Sys-NL} defined in $[t_0, t_1]$. Assume that the linearized system around $(\bar x, \bar u)$ is controllable, i.e., the control system \eqref{Sys-NL-L} is controllable for the time interval $[t_0, t_1]$.  Then the control system \eqref{Sys-NL} is locally controllable along $(\bar x, \bar u)$. 
\end{theorem}

\begin{proof}  
We present two solutions.

\noindent 
{\it Solution 1}: For each $\delta x_0$ and $\delta x_1$ in a neighborhood of $0$ of $\mR^n$  let  $\delta u \in C^1([t_0, t_1])$ be defined by the construction given in the result of Kalman such that the $C^1$ trajectory $(\delta x, \delta u)$ on $[t_0, t_1]$ of the linearized system \eqref{Sys-NL-L} satisfying 
$$
\delta x (t_0) = \delta x_0 \quad \mbox{ and } \quad \delta x_1 (t_1) = \delta x_1. 
$$
The control $\delta u$ is explicitly given by, see \Cref{lem-Kalman-NA},  
\be \label{thm-LinearTest-u-**}
\delta u (s) = B(s)^\top R(t_1, s)^\top \mC^{-1} (\delta x_1 - R(t_1, t_0) \delta x_0) \mbox{ in } [t_0, t_1].
\ee

We thus have 
\be \label{thm-LinearTest-p1} 
\delta x ' = A(t) \delta x + B(t) \delta u \mbox{ in } [t_0, t_1] \quad \mbox{ and } \quad \delta x(t_0) = \delta x_0. 
\ee

Define 
$$
\cG(\delta x_0, \delta x_1) = (\delta x_0, \cF(\delta x_0, \delta x_1)): = \big(\delta x_0, x(t_1) - \bar x(t_1) \big), 
$$ 
where $x$ is the solution of the system 
$$
x' = f(x, \bar u + \delta u) \mbox{ in } [t_0, t_1],  \quad \mbox{ and } \quad x(t_0) = x_0 + \delta x_0, 
$$
and $\delta u$ is given by \eqref{thm-LinearTest-u-**}. 

Since $f$ is of class $C^1$, it follows that $\cF$ is  of class $C^1$ with respect to $\delta x_0$ and $\delta x_1$ in a neighborhood of $(0, 0)$.  We next compute $\cF_{\delta x_1}(0, 0)$. 

Consider $\delta x_0 = 0$ and write $x = \bar x + \tdeltax$. We have 
\be \label{thm-LinearTest-td1}
\tdeltax' = f_x (\bar x, \bar u) \tdeltax + f_u (\bar x, \bar u) \delta u + o(\tdeltax, \delta u) = A(t) \tdeltax + B(t) \delta u +  o(\tdeltax, \delta u) \mbox{ in } C^0([t_0, t_1]),    
\ee
and 
\be \label{thm-LinearTest-td2}
\tdeltax(t_0) = 0. 
\ee

Since $\delta x_0 = 0$, it follows that $\delta x$ is the solution of the system
\be \label{thm-LinearTest-d} 
\delta x ' = A(t) \delta x + B(t) \delta u \mbox{ in } [t_0, t_1] \quad \mbox{ and } \quad \delta x(t_0) = 0,  
\ee
where  
\be\label{thm-LinearTest-u} 
\delta u (s) = R(t_1, s) B(s) B(s)^\top R(t_1, s)^\top \mC^{-1} \delta x_1 \mbox{ in } [t_0, t_1]. 
\ee
From \eqref{thm-LinearTest-u}, we have 
\be
\|\delta u\|_{C^0([t_0, t_1])} \le C \| \delta x_1 \|.  
\ee
We derive from \eqref{thm-LinearTest-td1} and \eqref{thm-LinearTest-td2} that 
\be \label{thm-LinearTest-td-p1}
\| \tdeltax \|_{C^0([t_0, t_1])} \le C \| \delta x_1 \|. 
\ee
Combining \eqref{thm-LinearTest-td1}, \eqref{thm-LinearTest-td2}, and \eqref{thm-LinearTest-d}, and using \eqref{thm-LinearTest-td-p1}, we obtain 
\be \label{thm-LinearTest-td-d}
\| \tdeltax - \delta x\|_{C^0([t_0, t_1])} = o( \| \delta x_1 \|). 
\ee
This implies 
$$
\cF(0, \delta x_1) - \cF(0, 0) = \tdeltax (t_1) = \delta x (t_1) + o(\delta x_1). 
$$
It follows that $\cF_{\delta x_1}(0, 0) = I$. 

We thus derive that 
$$
\nabla \cG(0, 0) = \left(\begin{array}{cc}
I_{n \times n} & 0 \\
* & I_{n \times n}
\end{array}\right). 
$$
The conclusion is now a consequence of the inverse theorem after noting \eqref{thm-LinearTest-u}.

\medskip 
\noindent{\it Solution 2}: Consider the map 
\be \label{def-Lambda}
\begin{array}{cccc}
\Lambda\colon  &  \overline{B_{\delta} (x_1)} & \to & \overline{B_{\delta} (x_1)} \\[6pt]
& \varphi & \mapsto & \varphi -  H(\varphi) +  x_1,  
\end{array}
\ee
where $H(\varphi) = x(t_1)$ with $x$ being the solution of the system 
$$
x' = f(x, \bar u + \delta u) \mbox{ in } [t_0, t_1] \quad \mbox{ and } \quad x(t_0) = x_0, 
$$
where, with $\delta x_1  = \varphi - \bar x (t_1)$ and $\delta x_0  = x_0 - \bar x (t_0)$,  
$$
\delta u (s) =   B(s)^\top R(t_1, s)^\top \mC^{-1} \big(\delta x_1 - R(t_1, t_0) \delta x_0 \big) \mbox{ in } [t_0, t_1]. 
$$
It is clear that 
\be \label{thm-LinearTest-pp1}
\| \delta u \|_{C^0([t_0, t_1])} \le C \Big( \| \delta x_1 \|  + \| \delta x_0 \| \Big).
\ee

Let $\delta x$ be the solution of the system 
$$
\delta x' = A(t) \delta x + B(t) \delta u \quad \mbox{ and } \quad  \delta x(t_0) = \delta x_0. 
$$
We have, by the choice of $\delta u$, 
$$
\delta x(t_1) = \delta x_1. 
$$

Let $\xi$ and $\txi$ in $B_{\rho} (x_1)$ and let $x$ and $\tx$ be the corresponding solutions with the corresponding controls $u(t) =\bar u + \delta u$ and $\tu (t) = \bar u + \tdeltau$. Denote $\delta x$ and $\delta \tx$ correspondingly. 

Then  
$$
H(\varphi) - H (\tvarphi) = x(t_1) - \tx(t_1). 
$$  

We have 
\be \label{thm-LinearTest-pp-dxtx}
(\delta \tx - \delta x)' = A(t) (\delta \tx - \delta x)  + B(t) (\delta \tu - \delta u)  \quad \mbox{ and } \quad (\delta \tx - \delta x) (t_0) = 0,  
\ee
which yields 
\be \label{thm-LinearTest-pp-dxtx-p1}
\|\delta \tx - \delta x \|_{C^0([t_0, t_1])} \le C \| \delta \tx_1 - \delta x_1 \|. 
\ee

We also have 
\be \label{thm-LinearTest-pp-xtx}
(\tx - x)' = \hat A(t) (\tx - x) + \hat B (t) (\delta \tu - \delta u) + o(\tx - x, \delta \tu - \delta u) \mbox{ in } C^0([t_0, t_1]), \mbox{ and }  (\tx - x)(t_0) = 0, 
\ee
where 
$$
\hat A(t) = \partial_x f(x(t), u(t)) \quad \mbox{ and } \quad \hat B(t) = \partial_u f(x(t), u(t)) \mbox{ in } [t_0, t_1]. 
$$
Thus
\be \label{thm-LinearTest-pp-xtx-p1}
\|\tx  - x\|_{C^0([t_0, t_1])} \le C \| \delta \tx_1 - \delta x_1 \|.  
\ee

Since 
$$
\|\hat A - A\|_{C^0([t_0, t_1])} + \|\hat B - B\|_{C^0([t_0, t_1])} \mbox{ is small as $\delta $ is small}, 
$$
we derive from \eqref{thm-LinearTest-pp-dxtx},   \eqref{thm-LinearTest-pp-dxtx-p1}, \eqref{thm-LinearTest-pp-xtx}, and \eqref{thm-LinearTest-pp-xtx-p1} that 
$$
\|(\tx  - x) - (\delta \tx - \delta x)\|_{C^0([t_0, t_1])} \le \frac{1}{2} \| \delta \tx_1 - \delta x_1 \| 
$$
as $\delta$ is small. The conclusion follows. 
\end{proof}

Here is a useful consequence 

\begin{corollary} Let $(x_e, u_e)$ be an equilibrium point. Assume that the linearized system around $(x_e, u_e)$ is controllable. Then the system is locally controllable around $(x_e, u_e)$. 
\end{corollary}

\begin{remark} \rm There is a proof of \Cref{thm-LinearTest} which deal with the whole $L^2$-class of controls. The proof is based on the implicit theorem for maps from an infinite dimensional normed space into a finite dimensional space. The version presented in this note takes into account the information from the proof of the Kalman rank condition. This allows one to deal only with finite dimensional spaces. The second solution is in the spirit of the first one but different. It is based on the idea that if one knows how to suitably build a control which brings $x_0$ at time $t_0$ to $x_1$ at time $t_1$ for the linearized system, then one can do it for the nonlinear system via $\Lambda$. This map has appeared in the study of KdV equation in the setting for which the linearized system is not even controllable \cite{CC04} (see also \cite{CKN24,Ng-KdV25}). 
\end{remark}

\begin{remark} \rm When the linearized system is not controllable, an powerful tool to deal with the nonlinear system in the finite dimensional setting is the use of the Lie brackets, see, e.g., \cite{Sontag98,Coron07}. 
In the infinite-dimensional setting, this idea is not easy to be adapted. Nevertheless,  there are several methods to analyze the controllability of a system even if the linearized system is not controllable such as the return method and the power series expansion. The structure of the nonlinearity plays a role then. For more details on these aspects, the reader is referred to \cite{Coron07}. 
\end{remark}

\section{Linear quadratic optimal control in finite and infinite horizons} \label{sect-optimal}

In this section, we study the theory of linear quadratic optimal control in finite ($0 < T  < + \infty$) and infinite horizons ($T = + \infty$) for the control system 
\be
x' = A x + B u \mbox{ in } (0, T) \quad \mbox{ and } \quad x(0) = \xi  \in \mR^n, 
\ee
associated with the cost function given by 
\be
\int_0^\top \Big(\|C x\|^2 + \|u \|^2 \Big) \, ds. 
\ee
Here $A$ is an $n \times n$ matrix, $B$ is an $n \times p $ matrix, and $C$ is a $m \times n$ matrix and $\| \cdot \|$ denotes the Euclidean norm. 

It is more general and more convenient for later purposes in the study of infinite horizon to study a slightly more general setting in the case $0< T< + \infty$ as follows. Let $0 < T < + \infty$,  and let $P_0$  be an $n \times n$ non-negative and symmetric matrix, we consider the cost function 
\be \label{def-J-0T}
J_T(x, u) = \int_0^T \Big(\|C x (s)\|^2 + \|u (s) \|^2 \Big) \, ds + \langle P_0 x(T), x(T) \rangle, 
\ee
and study the problem, for a given $\xi  \in \mR^n$,  
\be\label{inf-J-0T}
 \inf_{u \in L^2((0, T); \mR^p)} J_{T} (x, u),  
\ee
where $x$ is the solution of the system 
\be \label{sys-C-0T}
x' = A x + Bu \mbox{ in } (0, T) \quad \mbox{ and } \quad x(0) = \xi. 
\ee

We also denote, for $0 \le t \le T$,  
\be \label{def-J-tT}
J_{t, T} (x, u) = \int_t^T \Big(\|C x (s)\|^2 + \|u (s) \|^2 \Big) \, ds + \langle P_0 x(T), x(T) \rangle,  
\ee
and consider the problem, for a given $ \xi \in \mR^n$,  
\be \label{inf-J-tT}
 \inf_{u \in L^2((t, T); \mR^p)} J_{t, T} (x, u),  
\ee
where $x$ is the solution of the system 
\be \label{sys-C-tT}
x' = A x + Bu \mbox{ in } (t, T) \quad \mbox{ and } \quad x(t) = \xi. 
\ee

Here and in what follows, we often do not explicitly mention the relation \eqref{sys-C-tT} when we talk about \eqref{inf-J-tT} for notational ease. 

\medskip 
We first deal with the existence and uniqueness of the optimal control. 

\begin{proposition}\label{pro-Opt-WP} Let $T > 0$. Given $\xi \in \mR^n$, there exists a unique $\tu \in L^2((0, T); \mR^p)$ such that 
\be
J_T(\tx, \tu) =  \inf_{u \in L^2((0, T); \mR^p)} J_{T} (x, u),  
\ee
where $x$ is the solution of the system 
\be \label{pro-Opt-WP-sys}
x' = A x + Bu \mbox{ in } (0, T) \quad \mbox{ and } \quad x(0) = \xi. 
\ee
Here $\tx$ is the unique solution of \eqref{pro-Opt-WP-sys} with $u = \tu$. 
\end{proposition}

\begin{remark} \rm
After the proof of \Cref{pro-Opt-WP}, we often mention that $(\tx, \tu)$ is the optimal solution to insist the fact that $\tx$ is the corresponding solution to $\tu$ where $\tu$ is the optimal control (solution). 
\end{remark}

\begin{proof} The proof is standard and based on the strict convexity of the cost function and can be proceeded as follows. Let $(u_n) \subset L^2((0, T); \mR^p)$ be a minimizing sequence, i.e., 
\be
\lim_{n \to + \infty} J_T(x_n, u_n) = \inf_{u \in L^2((0, T); \mR^p)} J_{T} (x, u), 
\ee
where $x$ is the solution of \eqref{pro-Opt-WP-sys}. Here $x_n$ is the solution corresponding to $u_n$. It follows that $(u_n)$ is a bounded sequence in $L^2((0, T); \mR^p)$. Without loss of generality, one might assume that 
\be \label{pro-Opt-WP-p0}
u_n \rightharpoonup \tu \mbox{ in } L^2((0, T); \mR^p) \mbox{ for some } \tu \in L^2((0, T); \mR^p). 
\ee
Let $\tx$ be the solution corresponding to $\tu$. We claim that 
\be \label{pro-Opt-WP-p1}
\mbox{$x_n$ converges to $\tx$ pointwise in $[0, T]$},  
\ee 
and 
\be \label{pro-Opt-WP-p2}
\mbox{$x_n$ converges to $\tx$ pointwise in $L^2((0, T); \mR^n)$}.  
\ee
Claim \eqref{pro-Opt-WP-p1} is just a consequence of \eqref{pro-Opt-WP-p0} and the fact that 
$$
x_n(t) = e^{t A } \xi + \int_0^t e^{(t-s) A} B u_n(s) \, ds \quad \mbox{ and } \quad \tx(t) = e^{t A } \xi + \int_0^t e^{(t-s) A} B \tu(s) \, ds. 
$$
Claim \eqref{pro-Opt-WP-p2} is just a consequence of \eqref{pro-Opt-WP-p1} using the fact that $(x_n)$ is bounded in $L^\infty((0, T); \mR^n)$. 

As consequence of \eqref{pro-Opt-WP-p1} and \eqref{pro-Opt-WP-p2}, we have 
$$
\lim_{n \to + \infty } \int_0^T \|C x_n (s)\|^2  \, ds + \langle P_0 x_n(T), x_n(T) \rangle = \int_0^T \|C \tx (s)\|^2  \, ds + \langle P_0 \tx (T), \tx (T) \rangle,
$$
and as a consequence of \eqref{pro-Opt-WP-p0}, we obtain 
$$
\liminf_{n \to + \infty } \int_0^T \|u_n (s)\|^2  \, ds \ge \int_0^T \|\tu  (s)\|^2  \, ds. 
$$
It follows that 
$$
J_T (\tx, \tu) \le \liminf_{n \to + \infty} J_T(x_n, u_n). 
$$
Thus $\tu$ is a solution of the optimal control problem. 

We next deal with the uniqueness. Let $\tu$ and $\hu$ in $L^2((0, T); \mR^p)$ be two optimal solutions and let $\tx$ and $\hx$ be the two corresponding solutions. Since $P_0$ is symmetric, computations yield  
\begin{multline}\label{pro-Opt-WP-p3}
\frac{1}{2} \Big( J_T(\tx, \tu) + J_T(\hx, \hu) \Big) - J_T \left(\frac{1}{2} (\tx + \hx), \frac{1}{2} (\tu + \hu) \right) \\[6pt]
=\frac{1}{4} \int_0^\top \Big(\|C (\tx - \hx)(s)\|^2 + \|(\tu - \hu) (s) \|^2 \Big) \, ds + \frac{1}{4} \langle P_0 (\tx- \hx)(T), (\tx - \hx)(T) \rangle. 
\end{multline}
Note that $\frac{1}{2} (\tx + \hx)$ is the solution of \eqref{pro-Opt-WP-sys} corresponding to the control $\frac{1}{2} (\tu + \hu) $ and thus 
$$
\frac{1}{2} \Big( J_T(\tx, \tu) + J_T(\hx, \hu) \Big) - J_T \left(\frac{1}{2} (\tx + \hx), \frac{1}{2} (\tu + \hu) \right) \le 0. 
$$
We then derive from \eqref{pro-Opt-WP-p3} that $\tu = \hu$. The uniqueness follows. 
\end{proof}


Since the control system is linear, one can then prove that 
\be \label{thm-Opt-L}
\mbox{$\tu$ and $\tx$ are linear functions of $\xi$. }
\ee
To this end, we just need to note that if $(\tx_1, \tu_1)$ and $(\tx_2, \tu_2)$ are two optimal solutions 
corresponding to the initial conditions $\xi_1$ and $\xi_2$, respectively, then 
$(\tx_1 + \tx_2,  \tu_1 + \tu_2)$ and $(\tx_1 - \tx_2,  \tu_1 - \tu_2)$ are the optimal solutions corresponding to the initial conditions $\xi_1 + \xi_2$ and $\xi_1 - \xi_2$. This claim can be proved as follows. We have 
\begin{multline*}
J_T (\tx_1 + \tx_2,  \tu_1 + \tu_2) + J_T (\tx_1 - \tx_2,  \tu_1 - \tu_2) \\[6pt]= \int_0^T \Big( \|C (\tx_1 + \tx_2)\|^2 + \| u_1 + u_2 \|^2 \Big) \, ds + \langle P_0 (\xi_1 + \xi_2)(T), (\xi_1 + \xi_2 )(T) \rangle \\[6pt]
+ \int_0^T \Big( \|C (\tx_1 - \tx_2)\|^2 + \| u_1 - u_2 \|^2 \Big) \, ds  + \langle P_0 (\xi_1 - \xi_2)(T), (\xi_1 - \xi_2 )(T) \rangle
\\[6pt]
= 2 \int_0^T \Big( \|C \tx_1\|^2 +  \| C\tx_2)\|^2 + \| u_1\|^2  + \|u_2 \|^2 \Big) \, ds + 2 \langle P_0 \xi_1(T), \xi_1(T) \rangle +2 \langle P_0 \xi_2(T), \xi_2(T) \rangle. 
\end{multline*}
Thus  
\be \label{thm-Opt-L-p1}
J_T (\tx_1 + \tx_2,  \tu_1 + \tu_2) + J_T (\tx_1 - \tx_2,  \tu_1 - \tu_2)
= 2 J_T(\tx_1, \tu_1) + 2 J_T(\tx_2, \tu_2). 
\ee

Similarly, if $(x_+, u_+)$ and $(x_-, u_-)$ be such that $u_+, u_- \in L^2((0, T); \mR^p)$ and 
$$
x_{\pm}' = A x_{\pm} + B u_{\pm} \mbox{ in } (0, T) \quad \mbox{ and } \quad x_{\pm}(0) = \xi_1 \pm \xi_2, 
$$
then 
\be\label{thm-Opt-L-p2}
J_T (x_+, u_+) + J_T (x_-, u_-) 
= 2 J_T \Big(\frac{1}{2}(x_+ + x_-), \frac{1}{2}(u_+ + u_-) \Big) + 2  J_T \Big(\frac{1}{2}(x_+ - x_-), \frac{1}{2}(u_+ - u_-) \Big).   
\ee
Combining \eqref{thm-Opt-L-p1} and \eqref{thm-Opt-L-p2} yields
$$
J_T (x_+, u_+) + J_T (x_-, u_-) \ge J_T (\tx_1 + \tx_2,  \tu_1 + \tu_2) + J_T (\tx_1 - \tx_2,  \tu_1 - \tu_2)
$$
and the claim and the linearity follow. 

From \Cref{pro-Opt-WP} and \eqref{thm-Opt-L}, after changing the initial time, we derive that, for $ \xi \in \mR^n$,   
$$
\inf_{u \in L^2((\tau, T); \mR^p)} J_{\tau, T} (x, u), 
$$
where 
$$
x' = A x + Bu \mbox{ in } (\tau, T) \quad \mbox{ and } \quad x(\tau) = \xi,  
$$
is a quadratic function of $\xi$. 

For $\tau \in [0, T]$, let $P_T(\tau)$ be an $n \times n$ symmetric matrix and defined by 
\be  \label{def-Pt}
\langle P_T(\tau) \xi, \xi \rangle =  \inf_{u \in L^2((\tau, T); \mR^p)} J_{\tau, T} (x, u),  
\ee
where 
$$
x' = A x + Bu \mbox{ in } (\tau, T) \quad \mbox{ and } \quad x(\tau) = \xi. 
$$
The definition above make sense for $\tau \in [0, T)$ and $P_T(T)(\xi, \xi)$ is understood as $\langle P_0 \xi, \xi \rangle$ for $\xi \in \mR^n$.  Thus $\langle P_{T}(\tau)\xi, \xi \rangle$ is the cost from time $\tau$ to time $T$ for the initial data $\xi \in \mR^n$ at time $\tau$.

\medskip 
Here is the main result of the linear quadratic optimal control problem in finite horizon. 

\begin{theorem} \label{thm-Opt-T} Let $T > 0$. Given $\xi \in \mR^n$,  let  $\tu \in L^2(0, T)$ be the optimal control of \eqref{inf-J-0T} and \eqref{sys-C-0T}, 
and let $\tx$ be the corresponding solution. We have 
\be \label{thm-Opt-T-sys-opt1}
\left\{
\begin{array}{cl}\tx' = A \tx + B \tu & \mbox{ in } (0, T),\\[6pt]
\tx(0) = \xi,
\end{array} \right.
\ee
and
\be  \label{thm-Opt-T-sys-opt1-*}
\tu = - B^\top \ty  \mbox{ in } (0, T) \mbox{ where }  \left\{
\begin{array}{cl} \ty' = - A^\top \ty - C^\top C \tx & \mbox{ in } (0, T), \\[6pt]
\ty(T) = P_0 \tx(T).  
\end{array} \right. 
\ee
Moreover, 
\be \label{thm-Opt-T-sys-opt2}
\ty (t) = P_T(t) \tx (t)  \mbox{ in } [0, T]. 
\ee
Consequently, 
\be\label{thm-Opt-T-sys-opt3}
\tu (t) = - B^\top P_T(t) \tx (t) \mbox{ in } [0, T]. 
\ee
\end{theorem}

\begin{remark} \rm The variables $\ty$ is called the adjoint state of the state $\tx$ of the linear quadratic optimal control problem \eqref{inf-J-0T} and \eqref{sys-C-0T}.  
\end{remark}

Before giving the proof of \Cref{thm-Opt-T}, we state and prove the following basic result. 

\begin{lemma}\label{lem-Opt-T} Let $T>0$. Given $\xi \in \mR^n$, consider the pair $(\tx, \tu)$ in which $\tu \in L^2((0, T); \mR^p)$ and $\tx$ is the solution of the system 
$$
\tx' = A \tx + B \tu \mbox{ in } (0, T) \quad \mbox{ and } \quad \tx(0) = \xi. 
$$
Then $(\tx, \tu)$ is the optimal solution of \eqref{inf-J-0T} and \eqref{sys-C-0T} if and only if 
\be 
\tu (r)= - B^\top \Big(e^{(T-r)A^\top} P_0 \tx(T) + \int_r^T e^{(t-r) A^\top} C^\top C \tx \, dt\Big) \mbox{ in } (0, T). 
\ee
\end{lemma}

\begin{proof}[Proof of \Cref{lem-Opt-T}] Let $u \in L^2((0, T); \mR^p)$ and let $x$ be the solution of the system. 
$$
x' = A x + B u \mbox{ in } (0, T) \quad \mbox{ and } \quad x(0) = \xi. 
$$
We have 
\begin{multline} \label{thm-Opt-T-Diff-Energy}
J_T(x, u) - J_T(\tx, \tu) = \int_0^\top |C (x - \tx)|^2 + |u - \tu|^2 + \langle P_0 (x(T) - \tx(T)), x(T) - \tx(T) \rangle  \\[6pt]
+ 2 \int_0^\top \langle C^\top C \tx,  (x - \tx) \rangle  +\langle  \tu, u - \tu \rangle \, dt  + 2 \langle P_0 \tx(T), x (T) - \tx(T) \rangle. 
\end{multline}
Thus $(\tx, \tu)$ is the optimal solution if and only if 
\be \label{thm-Opt-T-p1}
\int_0^T \Big(\langle C^\top C \tx,  Lv  \rangle  +\langle  \tu, v \rangle \Big) \, dt  +  \langle P_0 \tx(T), L v(T) \rangle = 0 
\mbox{ for all } v \in L^2((0, T); \mR^p), 
\ee
where, for each $v \in L^2((0, T); \mR^p)$, we denote $Lv$ the  
solution of the system 
\be \label{thm-Opt-T-def-Lv}
(Lv)' = A L v + B v \mbox{ in } (0, T) \mbox{ with } L v(0) = 0. 
\ee
From \eqref{thm-Opt-T-def-Lv}, we have,  for $v \in L^2((0, T); \mR^p)$,  
\be \label{thm-Opt-T-Lv}
Lv(t) = \int_0^t e^{(t-s) A} B v(s) \, ds. 
\ee
It follows from \eqref{thm-Opt-T-Lv} that 
\begin{multline*}
\int_0^T \langle C^\top C \tx, Lv \rangle \, dt  = \int_0^T \langle C^\top C \tx (t),  \int_0^t e^{(t-s) A} B v(s) \, ds \rangle \, dt \\[6pt]
=\int_0^T  \int_0^t \langle B^\top  e^{(t-s) A^\top} C^\top C \tx (t),   v(s) \, ds \rangle \, dt , 
\end{multline*}
which yields 
\be \label{thm-Opt-T-p2}
\int_0^T \langle C^\top C \tx, Lv \rangle  
=  \int_0^T \langle B^\top  \int_s^\top e^{(t-r) A^\top} C^\top C \tx (t) \, dt, v(s) \rangle \, ds.
\ee
From \eqref{thm-Opt-T-Lv}, we also have
$$
 \langle P_0 \tx(T), Lv(T) \rangle =  \langle P_0 \tx (T),  \int_0^T e^{(T-s) A} B v(s) \, ds \rangle, 
$$
which yields 
\be \label{thm-Opt-T-p3}
 \langle P_0 \tx(T), Lv(T) \rangle =  \int_0^\top \langle B^\top e^{(T-s)A^\top} P_0 \tx(T), v(s) \rangle \, ds. 
\ee

Combining \eqref{thm-Opt-T-p1}, \eqref{thm-Opt-T-p2}, and \eqref{thm-Opt-T-p3}, we derive that \eqref{thm-Opt-T-p1} holds if and only if 
\be
\tu (r)= - B^\top  \Big(e^{(T-r)A^\top} P_0 \tx(T) + \int_r^T e^{(t-r) A^\top} C^\top C \tx \, dt\Big) \mbox{ in } [0, T].  
\ee
The conclusion follows from \eqref{thm-Opt-T-Diff-Energy}. 
\end{proof}

We are ready to give 

\begin{proof}[Proof of \Cref{thm-Opt-T}] We first prove \eqref{thm-Opt-T-sys-opt1-*}. Let $(\tx, \tu)$ be the optimal solution. By \Cref{lem-Opt-T}
\be \label{thm-Opt-T-tu-p0}
\tu (r)= - B^\top  \Big(e^{(T-r)A^\top} P_0 \tx(T) + \int_r^T e^{(t-r) A^\top} C^\top C \tx \, dt\Big) \mbox{ in } (0, T). 
\ee
Thus 
\be \label{thm-Opt-T-tu}
\tu(r)= - B^\top \ty(r) \mbox{ in } (0, T),  
\ee
where $\ty$ is the solution of the system 
\be \label{thm-Opt-T-ty}
\ty'=  - A^\top \ty - C^\top C \tx \mbox{ in } [0, T] \quad \mbox{ and } \quad \ty (T) = P_0 \tx(T). 
\ee

We now prove that 
$$
\ty(t) = P_T(t) \tx (t) \mbox{ in } [0, T],   
$$
where $P_T(t)$ is defined in \eqref{def-Pt}. In what follows in this proof, we denote $P_T(t)$ by $P(t)$ for notational ease.

We have 
\begin{multline*}
\frac{d}{dt} \langle \tx, \ty \rangle = \langle \tx', \ty \rangle +  \langle \tx, \ty' \rangle =  \langle A \tx + B \tu, \ty \rangle + \langle \tx, - A^\top\ty - C^\top C \tx \rangle \\[6pt]
\mathop{=}^{\eqref{thm-Opt-T-tu}} \langle A \tx - B B^\top \ty, \ty \rangle + \langle \tx, - A^\top\ty - C^\top C \tx \rangle = - |B^\top\ty|^2 - |C \tx|^2. 
\end{multline*}
Integrating the above identity from $\tau$ to $T$, we obtain, for $ 0\le  \tau \le  T$, 
$$
\langle \tx(T), \ty (T) \rangle  - \langle \tx(\tau), \ty (\tau) \rangle =  - \int_{\tau}^{T } \big(|C \tx|^2 + |\tu|^2  \big) \, ds, 
$$
which yields 
$$
\langle \tx(\tau), \ty (\tau) \rangle =   \int_{\tau}^{T } \big(|C \tx|^2 +  |\tu|^2 \big) \, ds  + \langle P_0 \tx(T), \tx(T) \rangle. 
$$
From the definition of $P(\tau)$ in \eqref{def-Pt}, we thus obtain
\be \label{thm-Opt-T-tx-ty1}
\langle \tx(\tau), \ty (\tau) \rangle =   \langle P(\tau) \tx(\tau), \tx(\tau) \rangle \mbox{ for } 0 \le \tau \le T. 
\ee

In particular, 
\be \label{thm-Opt-T-tx-ty1-0}
\langle \tx(0), \ty (0) \rangle =   \langle P(0) \tx(0), \tx(0) \rangle. 
\ee

One can verify from the uniqueness of the optimal control that $\ty(0)$ is a linear function of $\tx(0)$.  It follows that there exists an $n \times n$ matrix $M_0$ such that $\ty(0) = M_0 \tx(0)$. We claim that $M_0$ is symmetric and thus obtain 
\be \label{thm-Opt-T-tx-ty2-0}
\ty(0) = M_0 \tx(0)
\mbox{ for some $n\times n$ symmetric matrix } M_0. 
\ee
 
Indeed, let $(\tu_1, \tx_1)$ and $(\tu_2, \tx_2)$ be two optimal solutions corresponding to the initial data $\xi_1$ and $\xi_2$ at the time $t = 0$, respectively. We have, for $t \in (0, T)$,  
\begin{multline*}
\frac{d}{dt}\langle \tx_1(t), \ty_2 (t) \rangle =  \langle \tx_1', \ty_2 \rangle +  \langle \tx_1, \ty_2' \rangle 
= \langle A \tx_1 - B B^\top \ty_1, \ty_2 \rangle + \langle \tx_1, - A^\top \ty_2 - C^\top C \tx_2 \rangle \\[6pt]= - \langle B^\top \ty_1, B^\top \ty_2 \rangle - \langle C \tx_1,  C \tx_2 \rangle. 
\end{multline*}
and similarly,   for $t \in (0, T)$, we obtain 
$$
\frac{d}{dt}\langle \tx_2(t), \ty_1 (t) \rangle = - \langle B^\top \ty_2, B^\top\ty_1 \rangle - \langle C \tx_2,  C \tx_1 \rangle. 
$$
This implies
\be \label{thm-Opt-T-m1}
\langle \tx_1(T), \ty_2 (T) \rangle - \langle \tx_1(0), \ty_2 (0) \rangle = - \int_0^T \Big( \langle B^\top \ty_1, B^\top \ty_2 \rangle + \langle C \tx_1,  C \tx_2 \rangle \Big) \, ds 
\ee
and 
\be \label{thm-Opt-T-m2}
\langle \tx_2(T), \ty_1 (T) \rangle - \langle \tx_2(0), \ty_1 (0) \rangle = - \int_0^T \Big( \langle B^\top\ty_2, B^\top\ty_1 \rangle + \langle C \tx_2,  C \tx_1 \rangle \Big) \, ds.  
\ee

Since  $P_0$ is symmetric and, by \eqref{thm-Opt-T-ty}, 
$$
\langle \tx_1(T), \ty_2 (T) \rangle = \langle \tx_1(T), P_0 \tx_2 (T) \rangle, \quad  \langle \tx_2(T), \ty_1 (T) \rangle = \langle \tx_2(T), P_0 \tx_1 (T) \rangle, 
$$ 
we derive from \eqref{thm-Opt-T-m1} and \eqref{thm-Opt-T-m2} that 
$$
\langle \tx_1(0), \ty_2 (0) \rangle = \langle \tx_2(0), \ty_1 (0) \rangle. 
$$
Thus 
$$
\langle \xi_1, M_0 \xi_2 \rangle = \langle \xi_2, M_0 \xi_1 \rangle. 
$$
Since $\xi_1, \xi_2$ can be arbitrary in $\mR^n$, we derive that 
$$
\mbox{$M_0$ is symmetric} 
$$
and claim \eqref{thm-Opt-T-tx-ty2-0} is proved. 

By changing the starting time, we have thus proved that,  for $\tau \in [0, T)$, 
\be \label{thm-Opt-T-tx-ty2}
\ty(\tau) = M(\tau) \tx(\tau)
\mbox{ for some $n \times n$ symmetric matrix } M(\tau).  
\ee

As a consequence of \eqref{thm-Opt-T-tx-ty1} and \eqref{thm-Opt-T-tx-ty2}, we derive that  $
M(\tau) = P(\tau)$ in $[0, T)$, which yields
\be \label{thm-Opt-T-tx-ty}
\ty(\tau) = P(\tau) \tx (\tau) \mbox{ in } [0, T]. 
\ee

It is clear that \eqref{thm-Opt-T-sys-opt3} is a direct consequence of \eqref{thm-Opt-T-sys-opt1-*} and \eqref{thm-Opt-T-sys-opt2}. 

The proof is complete. 
\end{proof}

As a consequence of \Cref{thm-Opt-T}, we obtain the Riccati differential equation for $P_T(t)$. 

\begin{corollary}\label{cor-Opt-T} Let $T > 0$ and $P_0$ be an $n \times n$ symmetric matrix. Then $P_T$ satisfies the following system 
\be\label{RDE-T}
P_T'  + P_T A  +  A^\top P_T - P_T B B^\top P_T     + C^\top C =  0  \mbox{ in } [0, T]
\quad \mbox{ and } \quad P_T(T) = P_0.  
\ee
\end{corollary}

\begin{remark} \rm Note that if $P$ is a solution of \eqref{RDE-T} then $P^\top$ is also a solution since $P_0$ is symmetric. It follows that a solution of \eqref{cor-Opt-T} is automatically symmetric. System  \eqref{RDE-T} is a system of nonlinear ordinary differential equations. 
From the theory of the ordinary differential equations, the existence of a solution in the interval $[0, T]$ given the fact $P_T(T) = P_0$ is not trivial. 
\end{remark}

\begin{proof} Fix $ \xi \in \mR^n$. Let $(\tx, \tu)$ be the optimal solution and let $\ty$ be the corresponding adjoint state. We have, by \eqref{thm-Opt-T-sys-opt1} and \eqref{thm-Opt-T-sys-opt1-*} of  \Cref{thm-Opt-T}, 
$$ 
\left\{
\begin{array}{cl}\tx' = A \tx - B B^\top \ty  & \mbox{ in } (0, T),\\[6pt]
\tx(0) = \xi,
\end{array} \right. \quad \mbox{ and } \quad 
\left\{
\begin{array}{cl} \ty' = - A^\top \ty - C^\top C \tx & \mbox{ in } (0, T), \\[6pt]
\ty(T) = P_0 \tx(T), 
\end{array} \right. 
$$
Using the fact that $\ty(t) = P_T(t) \tx (t)$ and the equation of $\ty$, we obtain 
$$
(P_T \tx)' = - A^\top P_T \tx - C^\top C \tx \mbox{ in } [0, T]. 
$$ 
Since, by the equation of $\tx$,  
$$
(P_T \tx)' = P_T' \tx + P_T \tx' = P_T' \tx + P_T (A \tx - B B^\top \ty) = P_T' \tx + P_T (A \tx - B B^\top P_T \tx) \mbox{ in } [0, T],   
$$
it follows that 
$$
(P_T'  + P_T A  - P_T B B^\top P_T) \tx  =  (-A^\top P_T -C^\top C) \tx \mbox{ in } [0, T].
$$
Since $\xi$ is arbitrary, one can choose $\tx(t)$ arbitrary for $0 \le t \le T$. It follows that 
$$
P_T'  + P_T A  - P_T B B^\top P_T   = - A^\top P_T -C^\top C, 
$$
which yields \eqref{RDE-T}. 
 \end{proof}

\begin{remark}\rm Here is another standard way often given in the literature to establish \eqref{thm-Opt-T-sys-opt2} in \Cref{thm-Opt-T} after obtaining \eqref{thm-Opt-T-sys-opt1-*}. Since the dynamic equation is linear and the cost function is quadratic, one can expect that there is a relation 
$$
\ty(t) = P(t) \tx (t) \mbox{ for } t \in [0, T]. 
$$
As in the proof of \Cref{cor-Opt-T}, using \eqref{thm-Opt-T-sys-opt1-*}, one derives that $P$ satisfies \eqref{RDE-T} and consequently $P$ is symmetric since $P$ is symmetric. One can expect then that $(\tx, \tu)$ is given by the system 
\be \label{rem-Opt-T-p1}
\tx' = A \tx - B B^\top P(t) \tx \mbox{ in } (0, T) \quad \mbox{ and } \quad \tx(0) = \xi, 
\ee
and 
\be
\tu = - B^\top P(t) \tx \mbox{ in } [0, T]. 
\ee
It is indeed the case.  Set
\be
\ty(t) = P(t) \tx(t) \mbox{ in } [0, T].  
\ee
where $\tx$ is given by \eqref{thm-Opt-T-p1}.  We have 
\begin{multline*}
\ty'(t) = P' \tx + P \tx' \mathop{=}^{\eqref{RDE-T}} -(P A + A^\top P- P B B^\top P + C^\top C) \tx + P ( A \tx - B B^\top P(t) \tx) \\[6pt]
=  - A^\top P \tx  - C^\top C \tx =  - A^\top \ty  - C^\top C \tx. 
\end{multline*}
This implies 
$$
\ty(r) = e^{(T-r)A^\top} P_0 \tx(T) + \int_r^T e^{(t-r) A^\top} C \tx \, dt \mbox{ for } r \in [0, T]
$$
and this in turn yields that 
\be 
\tu (r)= - B^\top  \Big(e^{(T-r)A^\top} P_0 \tx(T) + \int_r^T e^{(t-r) A^\top} C \tx \, dt\Big). 
\ee
The conclusion now follows from \Cref{lem-Opt-T}. 
\end{remark}

Here is a useful consequence of \Cref{thm-Opt-T} dealing with infinite horizon. 

\begin{theorem} \label{thm-P} Assume that the finite cost condition holds for the triple $(A, B, C)$, i.e., for every $\xi \in \mR^n$, there exists $u \in L^2((0, + \infty); \mR^p)$ such that 
\be\label{pro-P-finite-cost}
\int_0^\infty \Big(\|C  x\|^2 + \|u \|^2 \Big) \, ds < + \infty, 
\mbox{ where } 
x' = Ax + Bu \quad \mbox{ and } \quad x(0) = \xi. 
\ee
Let $P$ be an $n \times n$ symmetric matrix  defined by \footnote{The existence of $P$ is a part of the conclusion.}
\be \label{pro-P-P}
\langle P \xi, \xi \rangle = \inf_{u \in L^2((0, + \infty); \mR^p)}  \int_0^\infty \Big(\|C  x\|^2 + \|u \|^2 \Big) \, ds
\mbox{ where } 
x' = Ax + Bu \mbox{  and  }  x(0) = \xi. 
\ee
Let $\tu$ be the optimal solution and $\tx$ be the corresponding solution. 
Then 
\be \label{thm-P-ux}
\tu = - B^\top P \tx \mbox{ in } [0, + \infty) 
\ee
and $P$ satisfies the algebraic Riccati equation 
\be \label{thm-P-P}
A^\top P + P A  - P B B^\top P + C^\top C = 0.  
\ee
\end{theorem}

\begin{remark} \rm The existence and uniqueness of the optimal solution in the infinite horizon is established as in the one in the finite horizon in \Cref{pro-Opt-WP}. The details are left to the reader. 
\end{remark}

\begin{proof} We first admit the existence of $P$ defined by \eqref{pro-P-P},  and establish \eqref{thm-P-ux} and \eqref{thm-P-P}.

Fix $T>0$ and consider the optimal control problem \eqref{inf-J-0T} and \eqref{sys-C-0T} with $P_0 = P$. Define $P_T(
\tau)$ for $0 \le \tau \le T$ by \eqref{def-Pt}. Then, by the choice of $P_0$ and the definition of $P$, one can show that  
\be \label{thm-P-observation1}
P_T(\tau) = P \mbox{ for } 0 \le \tau \le T. 
\ee
Given $\xi \in \mR^n$, let $(\tx_T, \tu_T)$ be the optimal solution of \eqref{inf-J-0T} and \eqref{sys-C-0T} with $P_0 = P$, and let $(\tx, \tu)$ be the optimal solution of \eqref{pro-P-finite-cost}. One can show that 
\be  \label{thm-P-observation2}
\tx_T = \tx \quad \mbox{ and } \quad \tu_T = \tu \mbox{ in } [0, T]. 
\ee
Applying \Cref{thm-Opt-T}, we have 
\be\label{thm-P-observation3}
\tu_T(t) = - B^\top P_T(\tau) \tx_T(t) \mbox{ in } [0, T]. 
\ee
Combining \eqref{thm-P-observation1}, \eqref{thm-P-observation2}, and \eqref{thm-P-observation3} yields 
\be \label{thm-P-p1}
\tu(t) = - B^\top P \tx_T(t) \mbox{ for } t \ge 0, 
\ee
since $T$ can be chosen arbitrary. This is \eqref{thm-P-ux}. 

Assertion \eqref{thm-P-P} is just a consequence of \eqref{cor-Opt-T} and \eqref{thm-P-observation1}.

It remains to prove the existence of $P$ and the proof is as follows. Given $T>0$, let $P_T(\tau)$ for $0 \le \tau \le T$ be defined by \eqref{def-Pt} with $P_0 = 0$. It is clear that 
$$
\langle P_T(0) \xi, \xi \rangle \mbox{ is increasing with respect to $T$.}
$$
Moreover, the family $\Big(\langle P_T (0) \xi, \xi \rangle
\Big)_{T >0}$ is bounded above by the finite cost condition. Therefore, $\langle P_T (0) \xi, \xi \rangle $ converges as $T \to + \infty$ for every $\xi \in \mR^n$. Since $P_T(0)$ is symmetric, it follows that 
$$
\langle P_T \xi, \eta \rangle = \frac{1}{4} \left( \langle P_T (\xi+\eta), (\xi+\eta) \rangle -  \langle P_T (\xi - \eta), (\xi - \eta) \rangle \right) \mbox{ for } \xi, \eta \in \mR^n. 
$$
Hence $\langle P_T (0) \xi, \eta \rangle $ converges as $T \to + \infty$ for every $\xi, \eta \in \mR^n$. Denote $a(\xi, \eta)$ the limit. It is clear that $a(x, y)$ is symmetric and bilinear. We next establish that $a$ is continuous. To this end, we first apply the Banach Steinhaus theorem \footnote{The proof of the Banach Steinhauss in the finite dimensional case is quite trivial!} to the family $P_T (0) \xi$ to obtain 
$$
\sup_{T > 0} \|P_T (0) \xi\|\mbox{ is finite for each $\xi$ in $\mR^n$.} 
$$
Applying the Banach Steinhaus theorem again to $P_T(0)$, we derive that 
$$
\sup_{T > 0} \|P_T \|_{\cL(\mR^n)} \mbox{ is finite}. 
$$
Hence $a$ is continuous. 

We thus just prove that there exists a bilinear continuous symmetric form $a$ defined on $
\mR^n \times \mR^n$ such that 
$$
\lim_{T \to + \infty} \langle P_T (0) \xi, \eta \rangle  = a (\xi, \eta) \mbox{ for all } \xi, \eta \in \mR^n.  
$$
Let $P$ be an $n \times n$ symmetric matrix such that 
$$
a(x, y) = \langle P x, y \rangle. 
$$
Then 
$$
\lim_{T \to + \infty} \langle P_T (0) \xi, \xi \rangle = \langle P \xi, \xi \rangle. 
$$
One can check that, for $ \xi \in \mR^n$, \footnote{Indeed, one can easily check that the LHS $\le $ the RHS. To check the fact the RHS $\le $ the LHS, one just needs to use the optimal solution of the RHS.}
$$
\lim_{T \to + \infty} \langle P_T (0) \xi, \xi \rangle = \inf_{u \in L^2((0, + \infty); \mR^p)}  \int_0^\infty \Big(\|C  x\|^2 + \|u \|^2 \Big) \, ds, 
$$
where $x$ is the solution of the system 
$$
x' = Ax + Bu \mbox{ in } [0, + \infty) \quad \mbox{ and } \quad x(0) = \xi. 
$$
The conclusion follows. 
\end{proof}

Here is an interesting application of the above results. 

\begin{proposition} The finite cost condition holds if there exists a feedback control $F$, i.e., $u = Fx$ such that the corresponding system 
$$
x' = (A + B F) x
$$
is exponentially stable. 
\end{proposition}

\begin{proof} Assume that the finite cost condition holds and $C = I$. Set 
$$
F = - B^\top P,  
$$
where $P$ is define in \Cref{thm-P}. Then the corresponding system is exponentially stable (see \Cref{sect-stability}) since
$$
\int_0^\infty |e^{t(A + BF)} \xi|^2 \, dt \mbox{ is finite for each $\xi \in \mR^n$.}
$$

It is clear that if there exists a feedback control $F$, i.e., $u = Fx$ such that the corresponding system 
$$
x' = (A + B F) x
$$
is exponentially stable, then the finite cost condition holds.  
\end{proof}

\begin{remark} \rm
System \eqref{thm-Opt-T-sys-opt1-*} in \Cref{thm-Opt-T} can be viewed as the Pontryagin maximum principle \cite{pontryagin2018mathematical} and \eqref{thm-Opt-T-sys-opt3} is in the spirit of dynamic programming due to Kalman \cite{bellman1966dynamic}. The analysis presented here combines these two important tools. This viewpoint seems not to be emphasized in the literature. 
\end{remark}

\begin{remark} \rm The approach given here give a quick way to derive the results in the infinite horizon \Cref{{thm-P}} from the ones in the finite results with the final cost directly \Cref{thm-Opt-T}. One approach frequently used in the literature is based on the existence and uniqueness of a positive symmetric solution of 
the algebraic Riccati equation \eqref{thm-P-P}. In the optimal control problem, the cost function $P_T$ or $P$ is more intrinsic than (solutions of) the Riccati equations \eqref{RDE-T} or \eqref{thm-P-P}. The viewpoint presented here is useful when one considers the infinite dimensional setting with unbounded control for which the theory on the Riccati equations are delicate and not completely complete, see, e.g., \cite{lasiecka1991differential} and the references there in. 
\end{remark}

\begin{remark} \rm Many parts of the results and approaches given in this section can be easily generalized to the infinite dimensional setting, see \cite{Ng-Trelat}. The materials given here have their roots from there. 
\end{remark}

\begin{remark} \rm The results given in this section can be found, see, e.g., \cite{lions1971optimal,zuazua2007controllability,Zabczyk08,TW09} and the references therein. 
\end{remark}

\section{Stabilization via Gramian operators} \label{sect-Gramian}

Consider the linear control system 
$$
x' = A x + B u.
$$ 
Assume that the control system is controllable. Then there is a  beautiful way to stabilize the system using the Gramian matrix, defined as follows, for $\lambda > 0$ sufficiently large,  
\be \label{def-Gramian}
Q = \int_0^\infty e^{- 2 \lambda t }e^{-t A} B B^\top e^{-t A^\top}  \, dt = \int_0^\infty e^{-t (A + \lambda I)} B B^\top e^{-t (A+\lambda I)^\top}  \, dt.  
\ee
We have, for $\xi_1, \xi_2 \in \mR^n$,  
$$
\langle Q \xi_1, \xi_2 \rangle = \int_0^\infty \langle B^\top e^{-t (A+\lambda I)^\top} \xi_1, B^\top e^{-t (A+\lambda I)^\top} \xi_2 \rangle  \, dt. 
$$
Note that $Q$ is invertible by the Kalman criterion or its proof. We next derive an interesting property of $Q$. We have, for $\xi_1, \xi_2 \in \mR^n$,  
\begin{multline*}
\langle Q \xi_1, (A + \lambda I)^\top \xi_2 \rangle + \langle (A + \lambda I)^\top  \xi_1, Q \xi_2 \rangle 
\\[6pt]
= \int_0^\infty \langle B^\top e^{-t (A+\lambda I)^\top} \xi_1 , B^\top e^{-t (A+\lambda I)^\top}  (A + \lambda I)^\top \xi_2 \rangle  \, dt 
\\[6pt] + \int_0^\infty \langle B^\top e^{-t (A+\lambda I)^\top} (A+\lambda I)^\top \xi_1 , B^\top e^{-t (A+\lambda I)^\top}  \xi_2 \rangle \, dt \\[6pt]
= -  \int_0^\infty \frac{d}{dt} \langle B^\top e^{-t (A+\lambda I)^\top} \xi_1 , B^\top e^{-t (A+\lambda I)^\top}  \xi_2 \rangle \, dt =  \langle B^\top \xi_1, B^\top \xi_2 \rangle. 
\end{multline*}
Consequently, we obtain  
\be \label{Gramian-c1}
A Q + Q A^\top + 2 \lambda Q - B B^\top = 0. 
\ee
Set 
$$
P = Q^{-1}. 
$$
Multiplying the equation \eqref{Gramian-c1} on the left and on the right by $P$, we get
\be \label{P-Gramian}
P A + A^\top P + 2 \lambda P - P B B^\top P = 0. 
\ee

We now construct a feedback using $P$. Let $x$ be the solution of the system 
$$
x' = Ax + Bu \mbox{ where } u = - B^\top P x, 
$$
i.e., 
$$
x' = (A  - B B^\top P) x. 
$$
Consider the Lyapunov function 
\be
V (x) := \langle P x, x  \rangle. 
\ee
One can easily show that 
\begin{multline}
\frac{d}{dt} V(x) =  \langle P x', x \rangle +  \langle P x, x' \rangle = \langle (P A - P B B^\top P)x , x \rangle +  \langle P x, (A - B B^\top P) x \rangle \\[6pt] 
= - 2 \lambda \langle P x, x \rangle - \langle B^\top P x, B^\top P x \rangle. 
\end{multline}
We then derive that 
$$
V(x(t)) \le e^{-2 \lambda t} V(x(0)) \mbox{ for } t \ge 0.  
$$
This implies  
$$
\|x(t) \| \le C e^{-\lambda t} \| x(0)\|. 
$$

We have just proved the following result. 
\begin{theorem} Assume that $(A, B)$ is controllable. Define $Q$ by \eqref{def-Gramian} and set 
$$
P = Q^{-1}. 
$$
There exists $C > 0$ such that every solution $x$ of the system 
$$
x' = Ax + Bu \mbox{ where } u = - B^\top P x. 
$$
satisfies 
$$
\| x(t) \| \le C e^{-\lambda t} \| x(0)\|. 
$$
\end{theorem}

\begin{remark} \rm Other choices of $Q$ are possible, see \cite{Ng-Riccati} for a discussion and related works. 
\end{remark}

\begin{remark} \rm There is a connection of \eqref{P-Gramian} with \eqref{thm-P-P}. In fact, \eqref{P-Gramian} can be viewed as a special case of \eqref{thm-P-P} with $A$ replaced by $A + \lambda I$ and $C =0$. 
\end{remark}

\begin{remark} \rm The idea in this section can be extended to the infinite-dimensional case, see \cite{Ng-Riccati} and the references therein, and can be extended even in the nonlinear setting with unbounded controls in the strongly continuous group setting but quite recently \cite{Ng-S-KdV, Ng-S-Schrodinger}. 
\end{remark}

\begin{remark} \rm The use of Gramians to study the stabilization was initiated by  Lukes \cite{Lukes68} and Kleinman \cite{Kleinman70} in the finite-dimensional setting. 
\end{remark}

\begin{remark} \rm There are others way to stabilize a control system. One way to do it is to use the backstepping approach which is also a powerful tool to stabilize a partial differential equations in one dimensional space, see, e.g., \cite{Coron07,Krstic08}. 
\end{remark}

\end{document}